\documentclass{amsart}
\usepackage{amsfonts, color, amssymb, graphicx}
\usepackage{amsmath,amscd,enumerate,amsthm}
\usepackage[latin1]{inputenc}
\usepackage[all,cmtip]{xy}
\usepackage{pictexwd,dcpic}

\newcommand{\scal}[1]{\langle #1 \rangle}
\newcommand{\virg}[1]{``#1''}

\def\ol{\overline}
\def\ul{\underline}
\def\wt{\widetilde}

\def\lra{\longrightarrow}
\def\lmt{\longmapsto}
\def\ra{\rightarrow}
\def\Ra{\Rightarrow}
\def\In{\subseteq}

\def\iff{\Leftrightarrow}

\def\PP{\Bbb{P}}

\def\RR{\Bbb{R}}
\def\ZZ{\Bbb{Z}}

\def\mc{\mathcal}
\def\mf{\mathfrak}

\def\proof{\textbf{Proof. }}

\def\be{\begin{equation}}
\def\ee{\end{equation}}
\def\qed{\begin{flushright}$\square$\end{flushright}}

\newtheorem{prop}{Proposition}[section]
\newtheorem{teor}[prop]{Theorem}
\newtheorem*{theor}{Theorem} %For the main thm
\newtheorem{cor}[prop]{Corollary}
\newtheorem{lem}[prop]{Lemma}

\theoremstyle{definition}

\theoremstyle{definition}
\newtheorem{defin}[prop]{Definition}

\theoremstyle{remark}

\def\V{\mc{V}}
\def\B{\mc{B}}
\def\H{\mc{H}}
\def\E{\mc{E}}
\def\P{\mc{P}}
\def\leaf{\mc{L}}
\def\K{\mc{K}}
\def\fol{\mc{F}}

\def\vf{\mf{X}}
\def\eps{\varepsilon}
\def\e{\varepsilon}
\def\g{\mf{g}}

\def\Iso{\textrm{Iso}\,}

\def\SO{\textrm{SO}}
\def\SU{\textrm{SU}}

\def\Sym{\textrm{Sym}}

\def\st{\,\Big|\;}
\def\codim{\textrm{codim\,}}
\def\rk{\textrm{rk\,}}
\def\span{\textrm{span\,}}
\def\im{\textrm{Im\,}}
\def\nablal{\nabla^{\leaf}}

\def\A{\mf{A}}
\def\B{\mf{B_{\leaf}}}
\def\C{\mf{C}}
\def\3{\bar{X}_3}
\def\4{\bar{X}_4}
\def\parl{\,/\!\!/\,}

\def\a{\alpha}
\def\b{\beta}
\def\Om{\Omega}
\def\X{\mc{X}}
\def\Y{\mc{Y}}
\def\Z{\mc{Z}}
\def\p{\tilde{p}}
\def\xx{\tilde{x}}
\def\ggamma{\tilde{\gamma}}
\def\J{\tilde{J}}
\newcommand{\Ohm}[1]{\Omega_{\Sigma,\Sigma}(S^{#1})}
\def\cond{$(\dagger)$}

\usepackage{cite}

\numberwithin{equation}{section}

\author{Marco Radeschi}
\title{Low Dimensional Singular Riemannian Foliations on Spheres}
\begin{document}
\maketitle
\begin{abstract}
Singular Riemannian Foliations are particular types of foliations on Riemannian manifolds, in which leaves locally stay at a constant distance from each other. Singular Riemannian Foliations in round spheres play a special role, since they provide \virg{infinitesimal information} about Singular Riemannian Foliations in general. In this paper we show that Singular Riemannian Foliations in spheres, of dimension at most 3, are orbits of an isometric group action.
\end{abstract}

\bigskip

\bigskip

A partition $\fol$ of a Riemannian manifold $M$ into complete properly immersed submanifolds, called \emph{leaves}, is a \emph{transnormal system} if any geodesic starting perpendicular to one leaf stays perpedicular to all the leaves it meets. A transnormal system is a \emph{Singular Riemannian Foliation} (SRF for short) if there exists a family of compactly supported smooth vector fields whose span, at each point $p\in M$, coincides with the tangent space of the leaf through that point.

Singular Riemannian Foliations were defined by Molino ( see \cite{Mo}) as a way to generalize the foliations obtained by the orbit decompositions of isometric group actions. Other special kinds of Singular Riemannian Foliations include the decomposition of a manifold $M$ into the fibers of a Riemannian submersion $\pi:M\to B$, or the partition of a space form by the parallel submanifolds of a given isoparametric hypersurface.

When a foliation is given by the orbits of an isometric group action, the foliation is called \emph{homogeneous}. The main result of our paper is:
\begin{theor}[Main Theorem]\label{TEOR}
Let $(S^n, \fol)$ be a Singular Riemannian Foliation in a round sphere, where the dimension of the leaves is $\leq 3$. Then the foliation is homogeneous.
\end{theor}

Singular Riemannian Foliations in spheres can be seen as the infinitesimal information of a generic SRF $(M,\fol)$ at a point. In fact, given a point $p$ in $M$, one defines an \emph{infinitesimal foliation}  $\fol_p$ to be a new SRF on the unit normal space $\nu^1_pL_p$ of the leaf at $p$. With the example of isometric group actions in mind, this construction generalizes the slice representation of the isotropy group $G_p$ on the normal space of an orbit $\nu_p(G*p)$.
This fact motivates the necessity of understanding Singular Riemannian Foliations in round spheres.

When all the leaves of a Singular Riemannian Foliation have the same dimension, the foliation is called a \emph{Regular Riemannian Foliation}. Regular Riemannian Foliations in spheres have been almost completely classified. 

In fact, a series of theorems by Ghys \cite{Ghys}, Haefliger \cite{Hae} and Browder \cite{Bro} shows that for Regular Riemannian Foliations the dimension of the leaves can only be 1,3 or 7. If the dimension is 1 or 3, Grove and Gromoll proved in \cite{GG} that the foliation is homogeneous. If the dimension is 7 then the sphere has dimension 15, and Wilking proved in \cite{W2} that the only foliation coming from a Riemannian submersion is the one given by the Hopf fibration $S^7\ra S^{15}\ra S^8$. It is conjectured that there are no other $7$ dimensional regular Riemannian foliations on $S^{15}$.
\\

When the foliation is not regular, the conditions become much less restrictive, and we are still far from a complete classification. Even the particular case of foliations with leaves of codimension one, which is equivalent to the theory of isoparametric hypersurfaces, is a deep and interesting subject in its own right (see \cite{Cec} for an expositions of the subject). These foliations have been classified, except in two special case. In particular, there are non homogeneous examples.

Another class of Singular Riemannian Foliations on spheres that has been successfully studied is the class of \emph{polar foliations}, i.e. foliations which admit a totally geodesic submanifold (the \emph{section}) that meets every leaf perpendicularly. Codimension one foliations are of course a special case. A theorem of Thorbergsson states that every irreducible polar foliation of codimension bigger than 1 is equivalent to the orbit decomposition of an s-representation i.e. the isotropy representation of a symmetric space, and is thus homogeneous.

All known \virg{irreducible} examples of SRF in spheres are either homogeneous, or isoparametric, or the Hopf fibration $S^7\ra S^{15}\ra S^8$. It is an interesting question to know if there are any other SRF on spheres.

The paper is organized as follows: in Section \ref{intro} we recall the basic concepts of Singular Riemannian Foliations. In Section \ref{prel} we prove some preliminary results about the geometry and topology of the natural stratification that arises on such foliations. Here we also prove the Main Theorem in the case of 1-dimensional Singular Riemannian Foliations. Section \ref{2dim} is devoted to the proof of the Main Theorem for 2 dimensional foliations, and Section \ref{cases}- Section \ref{last} to the case of 3 dimensional foliations. This last proof is divided into 5 cases: we list these cases in section \ref{cases}, and then proceed to cover them in separate sections.

\section*{ Acknowledgements. }This paper is part of the author's PhD thesis. The author would like to thank his advisor Wolfgang Ziller for his encouragement and many useful comments and suggestions. The author also wishes to thank Alexander Lytchak, for many insightful conversations.

\bigskip

\section{Preliminaries}\label{intro}
\bigskip
In this section we review the definition and basic properties of Singular Riemannian Foliations. For a more detailed exposition we refer the reader to \cite{Mo} and \cite{GW}.
\subsection{Regular Riemannian Foliations}\label{RRF}
If $M$ is a Riemannian manifold with metric $g$, a foliation is said to be a \emph{Regular Riemannian Foliation} if one of the following equivalent conditions is satisfied:
\begin{itemize}
\item around each point $p\in M$ there is a neighborhood $U_p$, a Riemannian manifold $(\ol{U}_p,\ol{g})$ and a Riemannian submersion $\pi_p:U_p\ra \ol{U}_p$ such that the foliation given by the connected components of the leaves in $U_p$ consists of the fibers of $\pi_p$. We will call $\ol{U}$, together with the projection $\pi_p$, the \emph{local quotient}, and the open sets $U_p$ are called \emph{simple neighborhoods}.
\item Any geodesic starting perpendicular to a leaf, stays perpendicular to all the leaves it meets.
\end{itemize}

The most basic example of a Regular Riemannian Foliations is the partition of $M$ into the fibers of a Riemannnian submersion $\pi:M\ra B$.  Such a foliation is called \emph{simple}. Another example is the decomposition of the orbits of an isometric action of a connected Lie Group $G$ on $M$, provided that all the orbits have the same dimension. Such a foliation is called \emph{homogeneous foliation}.

\subsection{Notation} For a Regular Riemannian Foliation $(M,\fol)$ there is a splitting $$TM=\H\oplus \V,$$ where $\V$, the vertical bundle, is the bundle of tangent spaces of the leaves, and $\H$, the horizontal bundle, is the bundle of normal spaces. When dealing with vertical vectors, we will use lower cased letters $u, v, w\in \V_p$ and for horizontal vectors we will use the letters $x,y,z\in \H_p$. The same letters, capitalized, will be used for vector fields: $U,V,W\in \vf(\V)$, and $X,Y,Z\in \vf(\H)$. When dealing with projections of a vector $b$ or a vector field $E$ onto the two bundles, we will use the notation $b^v,\,b^h,\, E^v,\, E^h$.

\subsection{Basic vector fields}
For a Regular Riemannian Foliation one defines \emph{basic vector fields} locally, as horizontal vector fields satisfying the following (equivalent) properties:
\begin{itemize}
\item On a local quotient $\pi_p:U_p\ra \ol{U}_p$, $X|_{U_p}$ is $\pi_p$-related to some vector field on $\ol{U}_p$.
\item for every vertical vector field $V\in \vf(\V)$, the bracket $[X,V]\in\vf(\V)$ as well.
\end{itemize}
Our definition of basic vector field coincides with the definition in \cite{GW}. Given an open set $U$, denote by $\mf{B}_U$ the set of horizontal basic vector fields on $U$. This is an infinite dimensional vector space. If we restrict our attention to a single leaf $\leaf$, we define $$\mf{B}_{\leaf}:=\left\{X|_{\leaf}\st X\in \mf{B}_U\right\}.$$ Notice that for every horizontal vector $x\in \H_p$ one can define (locally around $p$) a basic vector field $X\in \mf{B}_{\leaf}$ such that $X_p=x$, and if any two basic vector fields $X^1,X^2$ agree at a point $p$, then $X^1(q)=X^2(q)$ for all $q\in \leaf_p$, in a neighborhood of $p$ in the leaf. In particular, $\mf{B}_{\leaf}$ is finite dimensional.

One can use basic vector fields to define the so-called \emph{Bott connection} on the normal bundle of a leaf, which is the connection for which the basic vector fields are parallel (this is why basic vector fields are sometimes referred to as \emph{Bott parallel} vector fields, see for example \cite{W1}). Since basic vector fields have constant norm along a leaf, the Bott connection is compatible with the metric on the normal bundle. In particular, the holonomy group of the Bott connection lies in $O(q)$, where $q$ is the codimension of the foliation. Moreover, the Bott connection is locally flat (since it admits parallel local vector fields), and therefore the holonomy is discrete and only depends on the fundamental group of the leaf. The holonomy of the Bott connection is also called the \emph{isotropy group} of $\leaf$, and in the case of regular homogeneous foliations it coincides with the usual definition of isotropy.

\subsection{Tensors}
Given a Regular Riemannian Foliation, one defines the following tensors.
\begin{itemize}
\item The \emph{$S$-tensor}
\begin{eqnarray*}
S:\H\times\V&\lra&\V\\
(x,u)&\lmt& S_xu:=-\nabla_u^vx
\end{eqnarray*}
as the \emph{shape operator} of the leaves. One has $S\equiv 0$ iff the leaves are totally geodesic.
\item The \emph{$A$-tensor} (or \emph{O'Neill's tensor})
\begin{eqnarray*}
A:\H\times\V&\lra&\V\\
(x,y)&\lmt& A_xy:=\nabla_x^vy
\end{eqnarray*}
which is the obstruction for the horizontal distribution to be integrable. Thus $A\equiv0$ iff there exist submanifolds of dimension $=\codim \fol$ that meet every leaf transversely and perpendicularly. 
\end{itemize}
From these fundamental tensors, together with the metric, other tensors can be created:
\begin{itemize}
\item The $B$-tensor, defined by the formula $\scal{B(u,v),x}=\scal{u,S_xv}$, or more explicitly by $B(u,v):=\nabla_u^hv$.
\item The $A^*$-tensor, defined as the adjoint of the $A$-tensor $\scal{A_x^*u,y}=\scal{u, A_xy}$, or by the formula $A_x^*u:=-\nabla_x^hu$. Sometimes it is better to use the notation $A^ux:=A_x^*u$. 
\end{itemize}
As is well known, the $S$ tensor is symmetric, i.e. $B(u,v)=B(v,u)$, and the $A$ tensor is skew symmetric, i.e. $A_xy=-A_yx$.
Many relations hold between these tensors, their covariant derivatives, and the curvature operator, see for example \cite{ON}, or \cite{GW} pg. 44.

\subsection{Singular Riemannian Foliations}
A \emph{Singular Riemannian Foliation} (SRF for short) is a decomposition of $M$ into disjoint connected, complete, immersed submanifolds (not necessarily of the same dimension) such that:
\begin{itemize}
\item[-] Every geodesic meeting one leaf perpendicularly, stays perpendicular to all the leaves it meets.
\item[-] The foliation is a singular foliation, i.e. around each point $p\in M$ one can find local smooth vector fields spanning the tangent spaces of the leaves.
\end{itemize}
It is conjectured that the first assumption already implies the second.
Typical examples of SRF are:
\begin{itemize}
\item[-] foliations obtained by the orbits of an isometric group action. As in the regular case, these are called \emph{homogeneous foliations}.
\item[-] foliations obtained by taking the closures of leaves of a Regular Riemannian Foliation $\fol$, denoted by $\ol{\fol}$.
\end{itemize}

The \emph{dimension} of a foliation, denoted by $\dim \fol$, is the maximal dimension of leaves. The set of points whose leaf has maximal dimension is open, dense and connected in $M$. Such a set is called the \emph{regular part} of the foliation, and is denoted $M_{reg}$. An important and basic fact is that if $(M,\fol)$ is a SRF, then the restricted foliation $(M_{reg},\fol|_{M_{reg}})$ is a Regular Riemannian Foliation as defined before. In particular, one can still talk about basic vector fields, and the $A$ and $S$ tensors on $M_{reg}$. Notice though that they do not extend to the singular part of the foliation.

\subsection{New SRF from old}\label{newSRF}
Often new Singular Riemannian Foliations arise from elementary ones through sequences of basic manipulations.
\begin{itemize}
\item If $(M,\fol)$ is a foliation, and $U\In M$ is an open set, then $(U,\fol|_U)$ denotes the foliation on $U$ given by the connected components of the leaves in $\fol$ restricted to $U$. In other words, if $p\in U$ then the leaf through $p$ in $\fol|_U$ is the connected component through $p$ of $\leaf_p\cap U$.
\item If $(M,\fol)$ is a foliation, and $N\In M$ is a submanifold such that every leaf intersecting $N$ is contained in it, we say that $N$ is a \emph{saturated submanifold} and denote the restricted foliation by $(N,\fol|_N)$.
\item If $(M_i,\fol_i)$, $i=1,2$ are two foliations, then $(M_1\times M_2,\fol_1\times \fol_2)$ is the foliation given by the product of leaves. In other words, $\leaf^{M_1\times M_2}_{(p,q)}=\leaf^{M_1}_p\times\leaf^{M_2}_q$. Such a foliation is called \emph{product foliation}.
\item If $(S^{n_i},\fol_i)$, $i=1,2$ are two SRF on round spheres, one can construct the \emph{spherical join} $$(S^{n_1},\fol_1)\star (S^{n_2},\fol_2)=(S^{n_1+n_2+1},\fol_1\star \fol_2),$$ by setting $$\leaf^{\star}_{\cos(t)p+\sin(t)q}=\cos t\cdot \leaf^{S^{n_1}}_p+\sin t\cdot\leaf^{S^{n_2}}_q.$$
\item If $(S^n,\fol)$ is a SRF on a round sphere, one constructs a foliation $(\RR^{n+1},\fol^{hom})$ by setting $$\leaf^{hom}_v:=\left\{\begin{array}{ll} \|v\|\cdot \leaf_{{v\over\|v\|}}& v\neq 0\\ \{0\} & v=0\end{array}\right.$$
This is a SRF on $\RR^{n+1}$ for which $0$ is a closed leaf, and is called \emph{homothetic foliation} associated to $\fol$. Notice that by definition the homothetic foliation is invariant under homotheties centered at the origin. Conversely, if $(\RR^{n+1},\fol)$ is some SRF such that $0$ is a closed leaf, then every sphere centered at the origin is saturated, and if $\fol^1:=\fol|_{S^n}$, then $\fol=\left(\fol^1\right)^{hom}$.
\end{itemize}

\subsection{Stratification}
A SRF has a natural stratification. Given a positive integer $r$ we define the \emph{$r$-dimensional strata} of $(M,\fol)$ as the connected components of
\[
\left\{p\in M\st \dim\leaf_p=r\right\}=\bigcup_{\dim\leaf_i=r}\leaf_i.
\]
From now on we'll use the notation $\Sigma_r$ for a connected component of the $r$-dimensional stratum, and if $p$ is a point in $M$, $\Sigma^p$ will stand for the stratum containing $p$.
We now recall an important lemma due to Molino (cf. \cite{Mo}).
\begin{lem}[Homothetic transformation]\label{homtr}
Let $(M,\fol)$ be a SRF, $p\in M$ and $\Om_p\In \leaf_p$ a neighborhood of $p$ in the leaf through $p$. Let $\eps>0$ be such that the normal exponential map $\exp^{\perp}:\nu^{\eps}\Om_p\lra Tub_{\eps}(\Om_p)$ is a diffeomorphism onto the tubular neighborhood of $\Om_p$ of radius $\eps$. Then for each $\lambda\in(0,1)$ the homotetic transformation
\begin{eqnarray*}
h_{\lambda}:Tub_{\eps}(\Om_p)&\lra &Tub_{\eps}(\Om_p)\\
\exp_qv&\lmt&\exp_q(\lambda v)
\end{eqnarray*}
sends leaves to leaves.
\end{lem}
In particular it follows that
\begin{lem}
Every connected component of a stratum is a (possibly noncomplete) manifold. Moreover, every geodesic $\gamma(t)$ starting tangent to a stratum, and perpendicular to the leaf, stays in the stratum for $t\in (-\e,\e)$. Finally, the closure of a stratum $\Sigma_r$ is contained in the union of \virg{more singular strata:}
\[
\ol{\Sigma_r}\In \bigcup_{r'\leq r}\Sigma_{r'}
\]
\end{lem}
According to Lychak and Thorbergsson \cite{LT} more can be said. In fact, every geodesic starting tangent to a stratum and perpendicular to the leaf, stays in the stratum for all but a discrete set of times.

The stratum containing leaves of maximal dimension is just the regular part, $M_{reg}$. The stratum whose leaves have minimal dimension is called \emph{minimal stratum}. From the corollary above, it is clear that every component of the minimal stratum is a closed embedded manifold without boundary.

\subsection{Infinitesimal foliation} Given a SRF $(M, \fol)$ and a point $p\in M$, there exists a SRF $T_p\fol$ on $T_pM$ with the following properties:
\begin{itemize}
\item There is a neighborhood $O$ of $p$, and a diffeomorphism $\phi:O\lra T_pM$ onto the image, such that $\fol|_O$ is given by the preimages of of $T_p\fol$ under $\phi$.
\item $T_p\fol$ is preserved under homoteties: if $\lambda\in \RR\setminus 0$ and $v\in T_pM$ then $\leaf_{\lambda v}=\lambda \cdot \leaf_v$.
\item $T_p\fol$ only depends on the transverse metric. In other words if $\hat{g}$ is another metric such that $(M,\hat{g}, \fol)$ is again a SRF and the $\hat{g}$-distance between leaves is the same as the $g$-distance,  then $\widehat{T_p\fol}=T_p\fol$, where $\widehat{T_p\fol}$ denotes the foliation on $T_pM$ with respect to $\hat{g}$.
\end{itemize}
The SRF $(T_pM,T_p\fol)$ is called \emph{infinitesimal foliation} of $(M,\fol)$ at $p$. One can check that this foliation splits as
\[
(T_pM,T_p\fol)=(T_p\Sigma^p,\fol_1)\times(\nu_p\Sigma^p,\fol_2)
\]
where $\fol_1$ is just the foliation given by affine subspaces parallel to $T_p\leaf_p$, and $\fol_2$ has $0$ as a leaf. In particular $\fol_2={\fol\big|_{\nu_p^1\Sigma_p}}^{hom}$, where $\fol\big|_{\nu_p^1\Sigma_p}$ is the foliation in the unit sphere of $\nu_p\Sigma^p$ (see paragraph \ref{newSRF} above). This foliation $\fol\big|_{\nu_p^1\Sigma_p}$ contains all the information about $T_p\fol$ and will be called the essential foliation $\fol^{ess}_p$ at $p$. It might be worth remarking that $\fol^{ess}$ can also be defined in the following way: if $x\in \nu_p^1\Sigma^p$ then $$\leaf^{ess}_x:=\left\{y\in \nu_p^1\Sigma^p\st \exp_pty\in \leaf_{\exp_ptx}\quad \forall t\in (-\eps,\eps)\right\}.$$
By abuse of notation, we will always refer to $\fol^{ess}_p$ as the infinitesimal foliation, and denote it with $\fol_p$.

\subsection{Holonomy and Projectable Fields}
A SRF $(M,\fol)$ is defined by the existence of horizontal geodesics, i.e. geodesics that are perpendicular to the leaves they meet. In particular, by taking variations of horizontal geodesics around a fixed horizontal geodesic $\gamma$, one obtains special Jacobi fields which are called \emph{projectable Jacobi fields}. On a local quotient $\pi: U_p\ra \ol{U}_p$, they represent those Jacobi fields that project to Jacobi fields along $\pi(\gamma)\In \ol{U}_p$. On the regular part they are characterized by the formula
\[
\left(J'\right)^v=-A_{\gamma'}J^h-S_{\gamma'}J^v.
\]
Thus, in order to define a projectable Jacobi field, one only needs to specify $J(0)\in T_{\gamma(0)}M$ and $(J')^h\in \H_{\gamma(0)}$. Hence if $\P_{\gamma}$ denotes the vector space of projectable Jacobi fields along $\gamma$, then $\dim\P_{\gamma}=2\dim M-\dim \fol$.

A subset of $\P_{\gamma}$ is given by the so-called \emph{Holonomy Jacobi fields}, i.e. those Jacobi fields given by variations of geodesics which project to a fixed geodesic in the quotient. These holonomy Jacobi fields have important special properties:
\begin{itemize}
\item[-] They are always vertical.
\item[-] In the regular part, the space of holonomy Jacobi fields along a horizontal geodesic is a vector space, whose dimension is the maximal dimension of the leaves met by $\gamma$.
\item[-] In the regular part, the holonomy  Jacobi fields are characterized by the formula
\[
J'=-A_{\gamma'}^*J-S_{\gamma'}J.
\]
\item[-] On a local quotient, holonomy Jacobi fields are given by the kernel of $\pi_*:\P_{\gamma}\ra Jac_{\pi(\gamma)}$, where $Jac_{\pi(\gamma)}$ denotes the set of Jacobi fields along $\pi(\gamma)$, and $\pi_*(J)(t):=\pi_*\left(J(t)\right)$.
\item[-] Along $\gamma$, the vertical distribution $\V|_{\gamma}$ is spanned by the holonomy Jacobi fields.
\item[-] Suppose $\dim\leaf_{\gamma(t_1)<\dim\leaf_{\gamma(t_0)}}$, and consider $\Omega_{t_0}\In \leaf_{\gamma(t_0)}$ a neighborhood of $\gamma(t_0)$ in the leaf, as in the homothetic transformation Lemma \ref{homtr}. Consider the closest-point map $p:\Om(t_0)\ra \leaf_{\gamma(t_1)}$ and define $\Om(t_1):=p\left(\Om(t_0)\right)$. Then Molino \cite{Mo}, Lemma 6.1 states that the map $$p:\Om(t_0)\lra \Om(t_1)$$ is a submersion (non necessarily Riemannian). Moreover, the differential $p_*:T_{\gamma(t_0)}\Om(t_0)\ra T_{\gamma(t_1)}\Om(t_1)$ is given by $$p_*(v)=J_v(t_1),$$ where $J_v$ is the holonomy Jacobi field along $\gamma$ such that $J_v(t_0)=v$.
\end{itemize}
\subsection{SRF on spheres}\label{constancy}
In this section, we will recall some of the special properties that foliations in spheres share. Most of these properties still hold in more generic space forms.
\begin{prop}\label{const1}
Let $(S^n,\fol)$ be a SRF on a round sphere. Then:
\begin{itemize}
\item[-] Given two basic vector fields $X,Y$, the vertical vector $A_XY$ has constant norm along the leaves (i.e. the function $\|A_XY\|$ is basic).
\item[-] Given a basic vector field $X$, the shape operator $S_X$ has constant eigenvalues along a leaf.
\end{itemize}
\end{prop}
The first statement easily follows by O'Neill curvature formulas (\cite{ON}). The proof of the second statement, despite being well known, does not seem to appear in the literature. In \cite{GW}, a proof is given in the case of regular foliations. In the next section (see Proposition \ref{constpr}) we will provide a proof for the case of singular foliations.
In particular, if $x\in \H_p$ is a horizontal vector at a point $p\in S^n$ and $X$ is a basic vector field around $p$ with $X_p=x$, then the rank of $A_x:\H_p\ra \V_p$ is equal to the rank of $A_X:\mf{B}_{\leaf_p}\ra \vf(\leaf_p)$. In the same way, the eigenvalues of $S_x$ are the same as the eigenvalues of $S_X$.
\\

On spheres, and space forms in general, the O'Neill tensors satisfy very nice differential equations

\begin{prop}[\cite{GG}]\label{equations!}
Given a Riemannian foliation on a sphere, the following equations hold:
\begin{enumerate}
\item $(\nabla_X^vS)_X=S_X^2+Id -A_XA_X^*$.
\item $(\nabla_X^vA)_X=2S_XA_X$.
\end{enumerate}
\end{prop}

The following gives an important characterization of homogeneous foliations on spheres:
\begin{teor}[Homogeneity Theorem, \cite{GG}]\label{HomThm}
Let $(S^n,\fol)$ be a SRF on a round sphere. Suppose we can find a regular leaf $\leaf_0$, an open cover $\{U_{\a}\}$ of $\leaf_0$, and Lie subalgebras $\E_{\a}\In \vf(U_{\a})$, such that the following conditions hold:
\begin{itemize}
\item[a)] $\E_{\a}|_{U_{\a}\cap U_{\b}}=\E_{\beta}|_{U_{\a}\cap U_{\b}}$ whenever $U_{\a}\cap U_{\b}\neq\emptyset$,
\item[b)] the elements of $\E_{\alpha}$ span all of $T_p\leaf_0$, for all $p\in U_{\alpha}$,
\item[c)] If $V,W\in \E_{\a}$, then $\scal{V,W}$ is constant,
\item[d)] For every basic vector field $X\in \mf{B}_{\leaf_0}$, $S_X\E_{\alpha}\In \E_{\a}$ and $A^*_X\E_{\a}\In \E_{\a}$.
\end{itemize}
Then the SRF $(S^n,\fol)$ is homogeneous, i.e. there is a Lie group acting on $S^n$ by isometries, such that the leaves of $\fol$ correspond to the orbits of the group action. Moreover, this action is locally free on the regular part.
\end{teor}

In the case of a locally free group action $G\ra \Iso(S^n)$, $\E$ corresponds to the sheaf of invariant vector fields, i.e. local vector fields $V$ such that $g_*V=V$ for all $g\in G$ whenever well defined.

The theorem gives rise to the strategy we will adopt to prove our main theorem: we will consider several cases, and for each case we will prove homogeneity by creating a sheaf $\E$ that satisfies the conditions $(a)-(d)$.

In the same paper, conditions are given to find such Lie subalgebras.
\begin{teor}[\cite{GG}]\label{HomThm2}
Let $(S^n,\fol)$ be a SRF on a round sphere.
\begin{enumerate}
\item If there is a regular point $p\in S^n$, and a horizontal vector $x\in\H_p$ such that $A_x:\H_p\ra \V_p$ is surjective (or equivalently, $A_x^*$ is injective), then the sheaf
\[
\A:=\left\{A_XY\st X,Y\in \mf{B}_{\leaf_0}\right\}
\]
satisfies properies (a), (b), (d) of the Homogeneity Theorem.
\item If $\dim\fol\leq 3$, then the sheaf $\A$ defined above satisfies condition (c) of the Homogeneity Theorem.
\end{enumerate}
\end{teor}

If a SRF satisfies the condition of part 1 of the theorem above, the foliation is called \emph{substantial}. One of the main points of \cite{GG} is proving that a regular Riemannian foliation with $\dim\fol\leq 3$ is substantial. Since the Homogeneity Theorem \ref{HomThm} is not explicitly stated in \cite{GG}, and since this is the crucial tool we use, we indicate a proof (which uses the same strategy as in \cite{GG}).
\\

{\bf Proof of \ref{HomThm}. } Fix a $U_{\a}\In \leaf_0$, and consider the Lie algebra $\g=\E_{\a}$. Let $G$ be the 1-connected group with Lie algebra $\g$: then there is a local action of a neighborhood of the identity $U_G\In G$ on $U_{\a}$ (defined via the flows of the vector fields in $\E_{\a}$), and if we fix $p\in U_{\a}$ there is a map $\iota_p:U_G\ra U_{\a}$ defined by $\iota_p(g)=g\cdot p$. The differential of $\iota_p$ sends the right invariant vector fields (call $\g_R$ the Lie algebra of right invariant vector fields) to $\E_{\a}$. Consider now the Lie algebra $\g_L$ of left invariant vector fields, and call $\mc{K}_{p}$ the image of $\g_L$ under $(\iota_P)_*$. If $K\in \mc{K}_{p}$, then it is killing along the leaf, since for any $V_1,V_2\in \E_{\alpha}$
\[
\scal{\nabla_{V_1}K,V_2}+\scal{\nabla_{V_2}K,V_1}=\scal{\nabla_{K}V_1,V_2}+\scal{\nabla_{K}V_2,V_1}=X\scal{V_1,V_2}=0.
\]
Here we used the fact that $[\E_{\a},\mc{K}_p]=\iota_p[\g_R,\g_L]=0$, and condition (c) of the theorem, namely that $\scal{V_1,V_2}$ is constant along the leaf. Since the elements of $\E_{\a}$ span the whole tangent space of the leaf by condition (b), the equation above says that $K$ is killing along the leaf. In particular $\phi^t:=\phi_K^t$ is an isometry for small $t$. Now define $\Phi^t$ on $\nu (U_{\alpha})$ as $\Phi^t(x_p)=X_{\phi^t(p)}$, where $X$ is the basic vector field such that $X_p=x_p$. Condition $(d)$ implies that $(\phi^t,\Phi^t)$  respect the normal connection and the shape operator: namely, for every tangent vector $v$ and normal vector field $\eta$,
\[
\Phi^t(\nabla^{\perp}_v\eta)=\nabla^{\perp}_{\phi^t_*v}\Phi^t\eta\qquad \phi^t_*(S_{\eta}v)=S_{\Phi^t\eta}(\phi^t_*v)
\]
Using the Fundamental theorem of submanifold geometry, this is enough to prove that there is a global isometry $\psi^t:S^n\lra S^n$ such that $\psi^t\big|_{U_{\a}}=\phi^t$, $\psi^t_*\big|_{\nu(U_{\a}}=\Phi^t$. Take now $\ol{K}:={d\over dt}\big|_{t=0}\psi^t$. This is a global vector field, such that $\ol{K}|_{U_{\a}}=K$. By condition $(a)$, one can show that $\ol{K}$ is always tangent to $\leaf_0$, and by condition $(d)$ again it is tangent to all the other leaves as well. For each $K\in \mc{K}$ one gets one such global Killing vector field $\ol{K}$. Since they span the whole tangent space of $\leaf_0$, the homogeneous foliation obtained by the flows of these vector fields coincide with the original foliation $\fol$ in a neighborhood of $\leaf_0$, and therefore the two foliations agree everywhere.
\qed

%------------------ PRELIMINARIES ---------------------------------------------- PRELIMINARIES --------------------------------------------

\bigskip

\section{The geometry of the strata}\label{prel}

\bigskip

In this section we will prove some results on the structure of singular strata, which will be used throughout the rest of the paper. In particular, we will prove the following:
\begin{prop}
Let $(S^n,\fol)$ be a SRF. Then:
\begin{itemize}
\item $\Sigma_0$ is a totally geodesic subsphere $S^k$. Moreover, $(S^n,\fol)$ splits as a spherical join $$(S^n,\fol)=(S^k,\fol_0)\star (S^{n-k-1},\fol_1),$$ where $\fol_0$ consists of only $0$-dimensional leaves, and $\fol_1$ contains no $0$-dimensional leaves.
\item Suppose $\fol$ does not contain $0$-dimensional leaves. Then $\Sigma_1$ consists of a disjoint union of totally geodesic subspheres, orthogonal to each other.
\item Any singular stratum is a minimal submanifold.
\end{itemize}
\end{prop}

From now on, unless otherwise specified, we will be considering a SRF $(S^n,\fol)$ on a round sphere, and we will omit this in the lemmas and propositions from now on.
\\

Consider a singular leaf $\leaf$, and a horizontal geodesic $\gamma$ starting at $p\in \leaf$ and leaving that stratum. Then, for some small $t_0$ there is  a neighborhood $\Omega_{\gamma(t_0)}\In \leaf_{\gamma(t_0)}$ of $\gamma(t_0)$, a neighborhood  $\Omega_p\In \leaf_p$ of $p$, and a submersion
\[
\pi:\Omega_{\gamma(t_0)}\lra \Om_p
\]
with nontrivial kernel defined by the closest-point map. Let $v\in\ker\pi_*$. Call $\psi(t)=\gamma(t_0-t)$ (so that now $\psi(t)$ is starting at $\Omega_{\gamma(t_0)}$) and let $J(t)$ the unique holonomy Jacobi field along $\psi$ such that $J(0)=v$. Then
\[
\left\{\begin{array}{l} J(0)=v\\
J'(0)=-A^*_xv-S_xv\end{array}\right.
\]
where $x=\psi'(0)$.
Since we are on a sphere,
\[
J(t)=\cos(t) E_1(t)+\sin(t) E_2(t),\quad E_1'=E_2'=0,\; E_1(0)=v,\, E_2(0)=-A^*_xv-S_xv
\]
Now, the norm of $J$ is
\begin{eqnarray*}
\|J(t)\|^2&=&\cos^2t\|E_1(t)\|^2+\sin^2t\|E_2(t)\|^2+2\sin t\cos t\scal{E_1(t),E_2(t)}\\
&=&\cos^2t\left(\|A^*_xv+S_xv\|^2\tan^2t-2\scal{v,A^*_xv+S_xv}\tan t+ \|v\|^2\right)
\end{eqnarray*}
and it is immediate to check that $J$ goes to zero at $t_0$ iff $A_x^*v=0$ and $S_xv={1\over \tan t_0} v$. We just proved the following lemma:
\begin{lem}\label{usefullem}
Let $(S^n,\fol)$ be a SRF, $\leaf$ a singular leaf, $\gamma$ a horizontal geodesic starting at $p\in \leaf$, and $\pi: \Omega_{\gamma(t_0)}\lra \Om_p$ the local submersion. Then
\[
\ker d\pi_{\gamma(t_0)}=\left\{v\in \V\st A^*_xv=0,\,S_xv={1\over \tan t_0}v\right\}
\]
where $x=-\gamma'(t_0)$.
\end{lem}

This allow us to give a proof of the second part of Proposition \ref{const1}:
\begin{prop}\label{constpr}
Let $(S^n,\fol)$ be a SRF on a round sphere. Then for any basic vector field $X$, the shape operator $S_X$ has constant eigenvalues along a leaf.
\end{prop}
\proof Let $\fol$ be the SRF on a sphere, $p\in \leaf$ a regular point, $x\in \H_p$ a horizontal vector, $X$ a (locally defined) basic vector field, with $X_p=x$. Moreover, let $v\in \V_p$ be an eigenvector of $S_x$, say $S_xv=\lambda v$. Take $J$ the projectable Jacobi field along $\gamma(t):=\exp_p tx$ such that $J(0)=v$, $J'(0)=-S_xv$. It will have a zero somewhere, say at $t_0$, and there are 2 possibilities: either $\gamma(t_0)$ is a regular point, or it is singular. By the argument of \cite{GW}, Theorem 1.1, we know that if $\gamma(t_0)$ is a regular point, then the corresponding eigenvalue has constant multiplicity along $\leaf$. Suppose then, that $\gamma(t_0)$ is a singular point.
\\

Take $\p$ another point in $\leaf$, let $\xx:=X_{\p}$, and $\ggamma(t):=\exp_{\p}t\xx$. We want to show that $\dim E_{\lambda}(S_x)=\dim E_{\lambda}(S_{\xx})$. By the proposition above, there are spaces $K_p\In E_{\lambda}(S_x)$, $K_{\p}\In E_{\lambda}(S_{\xx})$, defined as $K_p=E_{\lambda}(S_x)\cap \ker A_x^*$, $K_{\p}= E_{\lambda}(S_{\xx})\cap \ker A_{\xx}^*$, such that $\dim K_p=\dim K_{\p}=\dim \leaf - \dim \leaf_{\gamma(t_0)}$. We then want to show that
\[
\dim E_{\lambda}(S_x)/K_p=\dim E_{\lambda}(S_{\xx})/K_{\p}
\]
Moreover, if we define $$\K:=\left\{J\st J \textrm{ is Jacobi field along $\gamma$ and } J(t_0)=0\right\},$$ it also follows from the proposition above that $ev_0(\K)=K_p$, where $ev_0$ is the evaluation at $t=0$. In the same way one can define $$\tilde{\K}=\left\{\J\st \J \textrm{ is Jacobi field along $\ggamma$ and } \J(t_0)=0\right\}$$ and again, $ev_0\tilde{\K}=K_{\p}$.
\\

Now, for every $[v]\in E_{\lambda}(S_x)/K_p$, pick a representative $v$ in $E_{\lambda}(S_x)$. Again take the projectable Jacobi field $J_v(t)$ with $J_v(0)=v,\,J_v'(0)=-S_xv$, and look at $J$ in an interval of the form $(t_0-\e,t_0)$. On this interval, the geodesic $\gamma$ is on the regular part, so we can look at the projection $\pi(\gamma)$ on a local quotient. We can look at the projected vector field $\pi_*J_v$, and notice that $\lim_{t\ra t_0} \|J_v(t)\|=0$. Notice, moreover, that  $\pi_*J_v$ does not depend on the choice of representative $v$ we started with. Now, on $\ggamma|_{(t_0-\e,t_0)}$ consider a projectable Jacobi field $\J_v$ that projects to $\pi_*J_v$, and such that $\J_v(t_0)=0$. Such a $\J_v$ is uniquely defined, up to a Jacobi field in $\tilde{\K}$. In particular, $\J_v(0)$ is an eigenvector of $S_{\xx}$ with eigenvalue $\lambda$, and it is well defined up to an element in $ev_0\tilde{\K}=K_{\p}$. Therefore, the map
\begin{eqnarray*}
E_{\lambda}(S_x)/K_p&\lra&E_{\lambda}(S_{\xx})/K_{\p}\\
{}[v]&\lmt&[\J_v(0)]
\end{eqnarray*}
is well defined, and has an inverse obtained by inverting the roles of $p$ and $\p$. Therefore the two spaces have the same dimension, which is what we wanted to prove.
\qed
\begin{prop}\label{usefulprop}
Let $\leaf_0$,\, $\leaf_1$ be singular leaves, $\gamma:[0,1]\ra S^n$ be geodesic that minimizes the distance between $\leaf_0$ and $\leaf_1$, and let $p_i=\gamma(i)\in \leaf_i$,  $i=0,1$. Let $\leaf_t:=\leaf_{\gamma(t)}$. 
\begin{enumerate}
\item All the leaves $\leaf_t$, $t\in(0,1)$ have the same dimension $d$, so that they belong to the same stratum $\Sigma$.%; moreover, the holonomy of these leaves in $\Sigma$ is constant.
\item If $d(\leaf_0,\leaf_1)<\pi$ then
\[
d\leq \dim\leaf_0+\dim\leaf_1.
\]
\item If $d= \dim\leaf_0+\dim\leaf_1$, then the following are true:
\begin{itemize}
\item There are orthogonal subspaces $V_0, V_1\In \RR^{n+1}$ such that $\leaf_i\in V_i\cap S^n$, $i=0,1$. Equivalently, $d(\leaf_0,\leaf_1)=\pi/2$.
\item  The local submersions $\Om_t\to \Om_0$, $\Om_t\to \Om_1$, where $\Om_t\In \leaf_t$, are Riemannian submersions.
\item $\leaf_t$ locally splits as $\leaf_0\times\leaf_1$, and the local submersions correspond to the projections onto the corresponding factor.
\item $\leaf_t$ can be seen as the set
\[
\left\{\gamma(t)\st \gamma\textrm{ minimizing geodesic from $\leaf_0$ to $\leaf_1$}\right\}.
\]
\end{itemize}
\end{enumerate}
\end{prop}

\subsection*{Proof of (1).}
Suppose that that $\gamma$ goes through leaves of higher dimension. By the homothetic transformation lemma, all the leaves in between have the same dimension, hence they belong to the same stratum $\Sigma$. Now, for every small $\e>0$, $\gamma|_{[\e,1-\e]}$ is still minimizing and contained in the stratum.  %The constancy of holonomy is given by a standard trick: consider the leaves $\leaf_{\gamma(t)}$, for $t\in (\e, 1-\e)$, and suppose that the holonomy of some $\leaf_{\gamma(t_0)}$ does not fix $\gamma'(t_0)$

\subsection*{Proof of (2).} Fix a $t\in (0,1)$ and set $x=\gamma'(t)$, $p_{t}:=\gamma(t)$, $\leaf_t:=\leaf_{p_t}$. We know that we can find open sets $\Om_t\In \leaf_t,\Om_0\In \leaf_0,\Om_1\In \leaf_1$ around $p_t,p_0,p_1$ respectively, such that
\[
\pi_i:\Om_t\lra \Om_i,\quad i=0,1
\]
 are submersions. For a point $q\in \Om_t$, denote with $S_i(q)$ the preimage of $\pi_i^{-1}\left(\pi_i(q)\right)$, and with $E_i$ the distribution tangent to the $S_i(q)$, for points $q$. Moreover, for every $q\in \Om_t$ let $\gamma_q$ be the (unique) geodesic connecting $\pi_0(q)$ with $\pi_1(q)$. It is easy to check that $\gamma_q$ is horizontal, $\gamma_q(t)=q$, and if we define $X_q:=\gamma_q'(t)$, then $X$ is the basic vector field along $\Om_t$ such that $X_{p_t}=x$. From Lemma \ref{usefullem} we know that $E_i$ coincide with the distributions $\ker A_X\cap E_{1/\tan(t-t_i)}(S_X)$. In particular $E_0$, $E_1$ have intersection $0$ at all points, hence
 \begin{eqnarray*}
 \dim\leaf_0+\dim\leaf_1&=&(\dim\leaf_t-\dim E_0)+(\dim\leaf_t-\dim E_1)\\
 &=&\dim\leaf_t+(\dim\leaf_t-\dim E_0-\dim E_1)\geq \dim\leaf_t
 \end{eqnarray*}

\subsection*{Proof of (3).} Now suppose equality holds: then $E_0\oplus E_1=T\leaf_t|_{\Om_t}$. In particular,
\begin{equation}\label{eqisopar}
A_X\equiv 0,\quad E_i=E_{1\over \tan(t-t_i)}(S_X),
\end{equation}
and the maps
\begin{eqnarray*}
\pi_0|_{S_1(p_t)}:S_1(p_t)&\lra&\Om_0\\
\pi_0|_{S_0(p_t)}:S_0(p_t)&\lra&\Om_1\\
\end{eqnarray*}
are diffeomorphisms. Not just that: since $E_0\perp E_1$, then $E_1$ is the horizontal bundle of $\pi_0$, and vice versa. Moreover, if $v\in (E_1)_q$, then $(\pi_0)_*v={\sin(t_1-t_0)\over \sin t_1} P^{-t_0}v$, where $P^t$ denotes parallel translation along $\gamma_q$. In particular, if $v,w\in E_1$ then $$\scal{(\pi_0)_*v,(\pi_0)_*w}=\left({\sin(t_1-t_0)\over \sin t_1}\right)^2\scal{v,w},$$
and therefore $\pi_0$ is a Riemannian submersion, up to a factor ${\sin(t_1-t_0)\over \sin t_1}$. In the same way, $\pi_1$ is a Riemannian submersion as well, up to a factor ${\sin(t_1-t_0)\over \sin t_0}$.
\\

Pick a point $q_0\in \Om_0$, and $q_t\in \Om_t$ such that $\pi_0(q_t)=q_0$, and call $l_0=d(q_0,q_t)$. Then for all $S_0(q_t)$, we have $d(S_0(q_t),q_0)=d(q_t,q_0)=l_0$, and for all the $q_1\in \pi_1\left(S^0(q)\right)$, we have $d(q_1,q_0)=d(\leaf_1,\leaf_0)$. But since $\Om_1=\pi_1\left(S^0(q_t)\right)$, then every $q_0\in\leaf_0$ has the same distance to every other point in $\leaf_1$, and the same holds for every point in $\leaf_1$. In other words: for every point $q_0\in \leaf_0$, $q_1\in \leaf_1$, we have
\[
d(q_0,q_1)=cost=d(\leaf_0,\leaf_1)
\]
and this can happen if and only if this distance is $\pi/2$, and the 2 leaves lay in two orthogonal totally geodesic spheres.
\\

Finally, equation (\ref{eqisopar}) says that the basic vector field $X$ with $X_{\gamma(t)}=x$ is parallel (w.r.t. the normal connection). This means that $X$ is an isoparametric section, with $S_X$ having 2 distinct eigenvalues corresponding to the singular strata (see \cite{CO} for a definition of isoparametric section). It is known (see for example \cite{CO}) that in this case $\leaf_t$ is locally a product of the two eigendistributions of $S_X$ corresponding to the two fibers of $\pi_0$ and $\pi_1$. These two fibers, in turn, are isometric to $\leaf_1$ and $\leaf_0$, up to the factors ${\sin(t_1-t_0)\over \sin t_0}$ and ${\sin(t_1-t_0)\over \sin t_1}$ we discussed above.
\qed

\begin{cor}
As before, let $\leaf_0$, $\leaf_1$  be singular leaves, $\gamma(t)$ a minimizing geodesic between the two leaves, and $\leaf_t=\leaf_{\gamma(t)}$. If $\dim \leaf_t=\dim \leaf_0+\dim \leaf_1$, then $\leaf_0$, $\leaf_1$ belong to two different strata, at distance $\pi/2$ from each other.
\end{cor}
\proof Let $\Sigma_0,\Sigma_1$ be the singular strata containing $\leaf_0, \leaf_1$, respectively. It is enough to show that, for each leaf $\leaf$ in $\Sigma_1$
\[
d(\leaf_0,\leaf)=\pi/2.
\]
Let $\sigma(s)$ be a horizontal geodesic starting from $p_1=\gamma(1)\in \leaf_1$, and tangent to $\Sigma_1$.
Now, suppose that for small enough $s>0$, $d(\leaf_0,\leaf_{\sigma(s)})<\pi/2$. Take a sequence $s_n\ra0$, and let $\gamma_n:[0,\pi/2]\ra S^n$ be a minimizing geodesic between $\leaf_0$ and $\leaf_{\sigma(s_n)}$, starting at $p_0=\gamma(0)$. Let $x_n=\gamma_n'(0)$, and consider a converging subsequence of the $x_n$'s, that we will sill call $x_n$. Now, let $x_0$ be the limit of the $x_n$'s, and define $\gamma_0(t)=\exp(tx_0)$.

We know that $l(\gamma_n)=d(\leaf_0,\leaf_{\sigma(s_n)})<\pi/2$, and therefore
\[
\dim \leaf_{\gamma_n(t)}<\dim \leaf_0+\dim \leaf_1, \qquad\forall t\in (0,\pi/2).
\]
Taking the limit we know that the dimension of a leaf can only drop, therefore $\dim\leaf_{\gamma_0(t)}<\dim \leaf_0+\dim \leaf_1$. This contradicts the fact that $\gamma_0$, as a geodesic from $\leaf_0$ to $\leaf_1$, in $(0,1)$ meets only leaves of dimension equal to $\dim \leaf_0+\dim \leaf_1$.
Therefore, $d(\leaf_0,\leaf)=\pi/2$ for all the leaves $\leaf$ in a neighborhood of $\leaf_1$. This implies that the set of leaves $\leaf\in \Sigma_1$ such that $d(\leaf_0,\leaf)=\pi/2$ is open. Sinse this condition is closed as well, we obtain the result.
\qed

As a first easy corollary, we now prove that $\Sigma_0$ is a totally geodesic sphere.
\begin{cor}
The stratum $\Sigma_0$ consists of a totally geodesic sphere $S^k\In S^n$, $k\geq0$.
\end{cor}
From the homothetic transformation lemma, it is easy to see that each connected component is a totally geodesic sphere. Suppose there are two connected components $\Sigma_1$, $\Sigma_2$. Then either $d(\Sigma_1,\Sigma_2)=\pi$ (i.e. they consist of one point each, and they are antipodal, therefore $\Sigma=S^0$), or the minimizing geodesic connecting them goes through leaves of dimension $d\leq0+0=0$. But this means that the whole geodesic is in the 0-dimensional stratum, contradicting the hypotheses that $\Sigma_1$, $\Sigma_2$ were connected components of the stratum.
\qed
\begin{cor}\label{split}
Every SRF $(S^n,\fol)$ can be obtained as a join
\[
(S^k,\fol_0)\star(S^{n-k-1},\fol_1)
\]
where $\fol_1$ is a SRF without 0-dimensional leaves, and $\fol_0$ has only 0-dimensional leaves.
\end{cor}
\proof Given a SRF $\fol$ on $S^n$, suppose the stratum $\Sigma_0$ is nonempty, hence is a totally geodesic sphere $S^k$. Consider the sphere $S^{n-k-1}$ at distance $\pi/2$ from $\Sigma_0$: it is a saturated submanifold, since leaves stay at constant distance from $\Sigma_0$. Using the homotetic transformation lemma, one can see that $\fol$ is actually the join $(S^{n-k-1},\fol|_{S^{n-k-1}})\star(\Sigma_0, \fol|_{\Sigma_0})$.
\qed
\begin{cor}
Every SRF $(S^n,\fol)$ with 1-dimensional regular leaves, is homogeneous.
\end{cor}
\proof As in Corollary \ref{split}, $(S^n,\fol)=(S^{k},\fol_0)\star(S^{n-k-1},\fol_1)$, where $\fol_1$ is a regular 1-dimensional foliation, hence given by a $\RR$-action $\rho:\RR\lra \SO(n-k)$.

The initial foliation $\fol$ is then given by the $\RR$-action $\rho':\RR\lra\SO(n+1)$, where $\rho'$ is the composition of $\rho$ with the standard embedding $\SO(n-k)\In \SO(n+1)$.
\qed
\subsection{Minimality of singular strata}
The goal of this section is to prove that the singular strata are minimal submanifolds. This is the statement of Proposition \ref{minimality}. First we need to prove the following:
\begin{lem}\label{min}
Let $(S^n,\fol)$ be a SRF, $p\in \Sigma_r$, and $(\nu_p^1\leaf_p,\fol_p)$ the infinitesimal foliation at $p$. If $v, w\in \nu_p^1\leaf_p$ belong to the same (infinitesimal) leaf, then $\kappa(v)=\kappa(w)$, where $\kappa$ is the mean curvature form.
\end{lem}
\proof. Let us first assume that $v$ is a regular point for the infinitesimal foliation at $p$, and let $\gamma(t)$ be the horizontal geodesic with initial values $(p,v)$. Also, let $\leaf_t=\leaf_{\gamma(t)}$. Let $e_1,\ldots, e_r$ be an orthonormal basis of $\leaf_p=\leaf_0$, and let $E_1(t),\ldots, E_r(t)$ their vertically parallel extensions along $\gamma$. Also, let $J_{r+1}(t),\ldots J_{r_0}(t)$ be the holonomy Jacobi fields along $\gamma$ that vanish at $p$ (we are implicitly defining $r_0:=\dim\leaf_t,\,t\in(0,\varepsilon)$). Since we are on a sphere, we know that $E_i:=J_i/\|J_i\|$, $i=r+1\ldots r_0$ are vertically parallel, and $E_1(t),\ldots,E_r(t),E_{r+1}(t),\ldots E_{r_0}(t)$ form an othonormal basis of $T_{\gamma(t)}\leaf_t$, for $t\in(0,\varepsilon)$. Also, for $i=r+1,\ldots r_0$,
\[
S_{\gamma'(t)}E_i(t)={1\over \tan t} E_i(t)
\]
Therefore, we have
\begin{eqnarray*}
\kappa(v)&=&\sum_{i=1}^r\scal{\nabla_{e_i}v,e_i}=\lim_{t\ra 0+} \sum_{i=1}^r\scal{\nabla_{E_i(t)}v,E_i(t)}\\
&=&\lim_{t\ra 0+} \sum_{i=1}^{r_0}\scal{\nabla_{E_i(t)}v,E_i(t)}- \sum_{i=r+1}^{r_0}\scal{\nabla_{E_i(t)}v,E_i(t)}\\
&=&\lim_{t\ra 0+}\kappa(v_t)-{r_0-r\over \tan t}
\end{eqnarray*}
where $v_t=\gamma'(t)$.
Now, if $v,w$ are regular and belong to the same infinitesimal leaf, then $v_t,w_t$ represent, for $t$ small enough, the same basic vector field at different points. in particular, $\kappa(v_t)=\kappa(w_t)$, and from the computation above, $\kappa(v)=\kappa(w)$. The results now holds everywhere by continuity, since the regular set is dense.
\qed
\begin{prop}\label{minimality}
The singular strata are minimal.
\end{prop}
\proof Consider a point $p$ in a singular stratum $\Sigma$, and the singular leaf $p\in \leaf\In \Sigma$. We want to show that $\Sigma$ is minimal at $p$. But first, we want to establish a relation between the mean curvature of $\Sigma$ and the mean curvature of $\leaf$.

Let $S, \ol{S}$ the shape operators of $\leaf$, $\Sigma$ respectively. Let $v_1,\ldots,v_r$ an o.n. basis for $T_p\leaf$, and $w_{r+1},\ldots,w_t$ an o.n. basis for $\nu_p\leaf\cap T_p\Sigma$. For a horizontal vector $x\in \nu_p\Sigma$, the mean curvature of $\Sigma$ is
\begin{eqnarray*}
\kappa_{\Sigma}(x)&=&\sum_{i=1}^r\scal{\ol{S}_xv_i,v_i}+\sum_{j=r+1}^t\scal{\ol{S}_xw_{r+j},w_{r+j}}\\
&\stackrel{(*)}{=}&\sum_{i=1}^r\scal{S_xv_i,v_i}=\kappa_{\leaf}(x)
\end{eqnarray*}
Where $(*)$ holds because on the one hand $\scal{\ol{S}_xv_i,v_i}=\scal{S_xv_i,v_i}$, and on the other hand
\[
\scal{\ol{S}_xw_i,w_i}=\scal{x,\nabla^{\nu\Sigma}_{\phantom{\nu\Sigma}w_i}w_i}=0
\]
since geodesics starting normal to a leaf and tangent to the stratum, stay in the stratum.
Hence we have that
\[
\kappa_{\Sigma}=\kappa_{\leaf}\big|_{\nu\Sigma}
\]
In terms of mean curvature vectors: $H_{\Sigma}=pr_{\nu\Sigma}(H_{\leaf})$.

But by Proposition \ref{min}, we claim that the normalized mean vector $n_0=H_{\leaf}/\|H_{\leaf}\|$ lies in an infinitesimal leaf of dimension zero in $\nu_p^1\leaf$ (i.e. the leaf through $n_0$ consists only of $n_0$). In fact, if another vector $n'\in \nu_p^1\leaf$ lies in the same leaf of $n_0$, then
\[
\|H_{\leaf}\|=\scal{H_{\leaf}, n_0}=\kappa_{\leaf}(n_0)=\kappa_{\leaf}(n')=\scal{H_{\leaf},n'}=\|H_{\leaf}\|\cos\theta
\]
Thus $\theta=0$ and so $n'=n_0$.
This proves that $n_0$ lies on a zero dimensional leaf. This is equivalent to saying that $n_0$ lies in $T^1_p\Sigma$, hence that $H_{\leaf}\in T_p\Sigma$, and therefore
\[
H_{\Sigma}=pr_{\nu\Sigma}H_{\leaf}=0.
\]
And thus $\Sigma$ is minimal at $p$.
\qed
\subsection{Connectedness}

In this section we prove the following

\begin{prop}\label{connectedness}
Suppose $\Sigma$ is a $k$-dimensional compact stratum of $S^n$, such that every horizontal geodesic starting at $\Sigma$ meets $\Sigma$ again for the first time at distance $\pi$ or larger. Then $(S^n,\Sigma)$ is $k$-connected. In particular, $\Sigma$ is homeomorphic to $S^k$.
\end{prop}
\proof Let $\Ohm{n}$ be the set of paths $\gamma:[0,1]\lra S^n$ starting and ending in $\Sigma$, and let
\[
E:\Ohm{n}\lra \RR
\]
be the energy function. Recall the following proposition (cf. \cite{PP}, pag.180):
\begin{prop}
Let $M$ be a complete Riemannian manifold and $A\In M$ a compact submanifold. If every geodesic in $\Omega_{A,A}(M)$ has index $\geq k$, then $A\In M$ is $k$-connected.
\end{prop}
According to the proposition above it is enough to show that every geodesic in $\Ohm{n}$ has index at least $k$. Let $\gamma$ be such a geodesic. By assumption, this geodesic has length which is a multiple of $\pi$. Also, if it is longer than $2\pi$, the parallel vector fields $V_i$ along $\gamma$ starting tangent to $\Sigma$ give rise to a $k$-dimensional subspace where the energy has negative second derivative, and therefore that the index of $\gamma$ is at least $k$. Assume then that $\gamma$ has length $\pi$, i.e. it connects a point $p\in \Sigma$ with it is antipodal point $a(p)=-p$. Remember that the antipodal map $a$ is an isometry that preserves the foliation. In particular it preserves the vertical and horizontal spaces, the $A$ tensor and the $S$ tensor, every strata, and their shape operators.

Now consider the $k$-dimensional space of vector fields along $\gamma$, consisting of the parallel vector fields starting tangent to $\Sigma$. Notice that in this case, if $V$ is such a vector, then $V(\pi)=-a_*(V(0))\in T_{-p}\Sigma$, so these vector fields represent tangent vectors in $\Ohm{n}$ at $\gamma$. Of course the proposition will be proved if we show that the second variation of the energy along all these vector fields is negative. We compute it:
\[
E''(V,V)=-\int_0^1\scal{R(V(t), \gamma'(t))\gamma'(t), V(t)}dt+\scal{S_{\gamma'(\pi)}V(\pi),V(\pi)}-\scal{S_{\gamma'(0)}V(0),V(0)}
\]
Now, first of all $\scal{R(V(t),\gamma'(t))\gamma'(t),V(t)}=\pi\|V(0)\|^2$. Secondly,
\begin{eqnarray*}
\scal{S_{\gamma'(\pi)}V(\pi),V(\pi)}&=&\scal{S_{a_*\gamma'(0)}\left(-a_*V(0)\right),-a_*V(0)}\\
&=&\scal{-a_*\left(S_{\gamma'(0)}V(0)\right),-a_*V(0)}\\
&=&\scal{S_{\gamma'(0)}V(0),V(0)}
\end{eqnarray*}
so the boundary terms in the second variation equation cancel out. Therefore
\[
E''(V,V)=-\pi\|V(0)\|^2<0
\]
and the proposition is proved.\qed

\subsection{Geometry of $\Sigma_1$}\label{geomSigma1}

The main goal of this section is to prove Proposition \ref{sigma1}, which states that $\Sigma_1$ consists of a disjoint union of totally geodesic spheres.
First of all, since $\Sigma_0$ is empty, $\Sigma_1$ is compact, and is minimal by Proposition \ref{minimality}. We will first prove that it is a ruled submanifold in $S^n$. 
For this recall
\begin{defin}[Ruled submanifolds]
A \emph{ruled submanifold} $M^k\In \ol{M}^n$ is a submanifold which is foliated by $(k-1)$-dimensional totally geodesic submanifolds, i.e. for every point $p\in M$ there is a submanifold of $M$ passing through $q$, totally geodesic in $\ol{M}$ and of codimension 1in $M$.
\end{defin}
Now, notice that for every point $p\in\Sigma_1\In S^n$, one can ``exp out'' from $p$ the space $\nu_p(\leaf_p)\cap T_p\Sigma_1$ and obtain a totally geodesic sphere, contained in $\Sigma_1$ of codimension 1, passing through $p$.
Moreover, we can say the following:
\begin{prop}\label{myprop}
Every connected component of $\Sigma_1$ is homeomorphic to a sphere.
\end{prop}
\proof We will use Proposition \ref{connectedness}. Therefore, we need to prove that any horizontal geodesic leaving a certain connected component of $\Sigma_1$ goes back to that component at distance at $\pi$ or larger.

Suppose that there is a horizontal geodesic $\gamma:[0,1]\ra S^n$ starting at some connected component of $\Sigma_1$, ending at $\Sigma_1$ and shorter than $\pi$. Call $\leaf_0:=\leaf_{\gamma(0)}$, and $\leaf_1:=\leaf_{\gamma(1)}$. Then by Lemma \ref{usefulprop} the leaves in between must have dimension 2, and $dist(\leaf_0,\leaf_1)=\pi/2$. Moreover, the two strata containing $\leaf_0$ and $\leaf_1$ have mutual distance $\pi/2$, and in particular are not the same.

Therefore, every geodesic leaving some component of $\Sigma_1$ can only meet the same component  at distance $\pi$ or $2\pi$, and possibly some other components at distance multiple of $\pi/2$). Thus Corollary \ref{connectedness} finishes the proof.
\qed

Summing up, every connected component of $\Sigma_1$ is a compact, minimal, ruled submanifold of $S^n$, which is homeomorphic to a sphere.

Minimal ruled submanifolds in space forms were studied extensively, and a classification can be found in \cite{BDJ}. Before stating the classification, we recall the definition of \emph{generalized helicoid}:
\begin{defin}
Let $K$ be a Killing vector field on $S^n$, and $A(t)$ the one-parameter group of isometries generated by the flow of $K$. Moreover, let $S^k\In S^n$ be a totally geodesic submanifold which is perpendicular to $K$, and suppose that there is an open set $U\in S^k\times \RR$ such that the map
\begin{eqnarray*}
X:S^k\times \RR&\lra&S^n\\
(p,t)&\lmt&A(t)\cdot p
\end{eqnarray*}
is regular in $U$.
Then the map $X\big|_U$ is called \emph{generalized helicoid}.
\end{defin}
The main theorem in \cite{BDJ} classifies minimal ruled submanifolds in space forms. In particular:
\begin{teor}
Let $\Sigma^k$ be a minimal ruled submanifold of $S^n$. Then there exists a generalized helicoid
\[
X:S^{k-1}\times \RR\lra S^n
\]
and an open set $U\In S^{k-1}\times \RR$ such that $X$ restricted to $U$ parametrizes $\Sigma$.
\end{teor}

In particular if $\Sigma$ is complete, then the map is globally defined $X:S^{k-1}\times \RR\to \Sigma$.
\begin{prop}\label{sigma1}
Each connected component of $\Sigma_1$ is a totally geodesic sphere.
\end{prop}
\proof Consider a component $C\In\Sigma_1$. By what we said so far, we know that $C$ is a generalized helicoid. In particular there is a killing field $K$ whose restriction to $C$ is tangent to $C$, and a totally geodesic $k-1$-dimensional sphere $S^{k-1}$ in $C$, whose orbit under the flow of $K$ generates all of $C$. Moreover, from \cite{BDJ}, Proposition 3.20, we know that there exists a closed geodesic $\gamma$ in $\Sigma_1$ which is preserved by $K$ (i.e. $\gamma$ is an integral curve of $K$). We will identify the closed geodesic $\gamma:[0,2\pi]\lra S^n$ with its image.
\\

For each $t\in [0,2\pi]$, let $\wt{S}_t$ the totally geodesic $(k-1)$-sphere in $C$, passing through $\gamma(t)$. By definition of a generalized helicoid, $K$ is everywhere orthogonal to $\wt{S}_t$.
Now, use the following notation:
\begin{itemize}
\item We write $S^n$ as a spherical join $S^n=\gamma\star S^{n-2}$.
\item We define $S_t:=\wt{S}_t\cap S^{n-2}$. $S_t$ is a $(k-2)$-sphere, and $\wt{S}_t=S_t\star \{\pm \gamma(t)\}$.
\item For each $p\in S_t\In S^{n-2}$, let $\psi_t:[0,\pi/2]\lra S^n$ be the unit speed minimizing geodesic from $p$ to $\gamma(t)$, and let $\vec{t}_p:=\psi_t'(0)$. This defines a parallel unit vector field $\vec{t}$ on $\nu(S^{n-2})|_{S_t}$.
\end{itemize}
The geodesic $\gamma$ can be chosen so that $K$ is tangent to the $S^{n-2}$ orthogonal to $\gamma$. Therefore, the flow of $K$ gives a one parameter family $t\lmt S_t$  of totally geodesic $(k-1)$-spheres in $S^{n-2}$. Notice that if $K$ is tangent to some $S_t$, then the spheres $S_t$ are all the same sphere $S$, and $C$ is a join $\gamma\star S$, which is a totally geodesic sphere.
Now  suppose that for some $t_1<t_2$, $t_2-t_1\neq \pi$, $S_{t_1}\cap S_{t_2}\neq 0$, and suppose $p\in S_{t_1}\cap S_{t_2}$. Then at that point, we have that $T_pC$ can be written as
\[
T_pC=\vec{t}_1\oplus \vec{t}_2\oplus T_pS_{t_1}=\vec{t}_1\oplus \vec{t}_2\oplus T_pS_{t_2}
\]
(we are using here that $t_2-t_1\neq\pi$, to know that $\vec{t_1}$ and $\vec{t_2}$ are linearly independent). But the term $\vec{t}_1\oplus \vec{t}_2$ belongs to $\nu(S^{n-2})$, so this is possible only if $T_pS_{t_1}=T_pS_{t_2}$, which means $S_{t_1}=S_{t_2}$. Moreover $K_p$ belongs to $T_pC\cap T_pS^{n-2}$, and by what we said before
\[
T_pC\cap T_pS^{n-2}=T_pS_{t_1}=T_pS_{t_2}.
\]
Therefore, $K_p\in T_pS_{t_1}$, and the same thing can be said about all the other points of $S_{t_1}$. It follows that $K|_{S_{t_1}}$ is tangent to $S_{t_1}$, and therefore $C$ is a totally geodesic sphere.
\\

The only possibility left is the case $S_{t_1}\cap S_{t_2}=\emptyset$ whenever $t_2-t_1\neq \pi$, which implies $\wt{S}_{t_1}\cap \wt{S}_{t_2}=\emptyset$. On the other hand, if $t_2-t_1=\pi$ it is easy to see that $\wt{S}_{t_1}=\wt{S}_{t_2}$. Let
\[
\sigma_t:\wt{S}_t\lra \wt{S}_{t+\pi}
\]
be the flow of $K$ at $t=\pi$. Then $C$ is diffeomorphic to
\[
S^{k-1}\times_{\ZZ} \RR
\]
where the generator of $\ZZ$ sends $(p,t)$ to $(\sigma_t(p),t+\pi)$. But this is a contradiction, since $C$ would not even be simply connected.
\qed
\subsection{Constancy properties of $A$- and $S$-tensors}
In this section, we slightly generalize lemma 4.3.2 of \cite{GW}.
\begin{prop}\label{Slem}
Suppose $M$ is a space form with a foliation $\fol$. Then, in the regular part:
\begin{enumerate}
\item $\scal{S_X^kA_XY,A_XZ}$ is a basic function for every basic $X$, $Y$, and every $k\in \ZZ_{\geq0}$.
\item If $S_X$ has all distinct eigenvalues, and eigenvectors $V_1,\ldots V_n$ then for any $Y\in \B$, $\scal{A_XY,V_i}$ is basic for all $i$. 
\end{enumerate}
\end{prop}
\subsection*{Proof of (1).}
First of all notice that in constant curvature $\scal{A_XY,A_XZ}$ is basic for every basic vector fields $X,Y,Z$. Moreover, notice that if $f$ is a basic function and $X$ is a basic vector field then $X(f)$ is again basic, since for every $V$ vertical vector,
\[
VX(f)=XV(f)+[X,V](f)=0+0=0.
\]
We now want to prove by induction that $\scal{S_X^kA_XY,A_XZ}$ is basic, or equivalently that $A_X^*S_X^kA_XY$ is a basic vector field.
Suppose we know it already up to $k$ (the case $k=0$ is true by what we said at the beginning). Then $\scal{S_X^kA_XY,A_XZ}$ is basic, and so is $X\scal{S_X^kA_XY,A_XZ}$. Assume now that $\nabla_XX=0$ at a leaf $\leaf$. We can then compute $X\scal{S_X^kA_XY,A_XZ}$ at $\leaf$:
\begin{eqnarray}
X\scal{S_X^kA_XY,A_XZ}&=&\scal{\nabla_X\left(S_X^kA_XY\right),A_XZ}+\scal{S_X^kA_XY,\nabla_X\left(A_XZ\right)}\\
&=&\sum_{i=0}^{k-1}\scal{S_X^iS'_XS_X^{k-i-1}A_XY,A_XZ}+\nonumber\\
&&+\scal{S_X^kA_X'Y,A_XZ}+\scal{S_X^kA_X(\nabla_X^hY),A_XZ}+\nonumber\\
&&+\scal{S_X^kA_XY,A'_XZ}+\scal{S_X^kA_XY,A_X\nabla_X^hZ}\nonumber
\end{eqnarray}
where the notation $T'$ was means $\nabla_X^vT$. In the equation before, the terms $\scal{S_X^kA_XY,A_X\nabla_X^hZ}$ and $\scal{S_X^kA_X(\nabla_X^hY),A_XZ}$ are basic because of the induction hypothesis and because $\nabla_X^hY$ is basic if $X,Y$ are. In particular the function
\begin{eqnarray}
\sum_{i=0}^{k-1}\scal{S_X^iS'_XS_X^{k-i-1}A_XY,A_XZ}+\scal{S_X^kA_X'Y,A_XZ}+\scal{S_X^kA_XY,A'_XZ}\label{eq}
\end{eqnarray}
is basic.
Also, remember that the following hold (\cite{GW}, pag. 44 and pag. 149):
\begin{eqnarray}
A'_XY&=&2S_XA_XY\label{eq1}\\
S'_XV&=&c\|X\|^2V+S_X^2V-A_XA_X^*V\label{eq2}
\end{eqnarray}
Substituting the equations (\ref{eq1}) and (\ref{eq2}) into equation (\ref{eq}), we get that
\begin{eqnarray*}
kc\scal{S^{k-1}_XA_XY,A_XZ}+(k+2)\scal{S_X^{k+1}A_XY,A_XZ}-\sum_{i=0}^{k-1}\scal{A_XA_X^*S_X^{k-i-1}A_XY,A_Z}
\end{eqnarray*}
is basic. The first term is basic by inductive hypothesis, and each term in the sum on the right is basic, since we can write the terms as
\[
\scal{A^*_XS_X^{k-i-1}A_XY,A_X^*A_XZ}
\]
that again by inductive hypothesis is the inner product of two basic vectors. In particular, the only remaining term
\[
(k+2)\scal{S_X^{k+1}A_XY,A_XZ}
\]
is basic, too.

\subsection*{Proof of (2).}
Write $A_XY$ as a linear combination of the eigenvectors of $S_X$
\[
A_XY=\sum_{i=1}^na_iV_i
\]
then
\[
S_X^kA_XY=\sum_{i=1}^{n}a_i\lambda_i^kV_i
\]
and if we look at a single leaf we get equations
\[
c_k=\scal{S_X^kA_XY,A_XY}=\sum_{i=1}^na^2_i\lambda_i^k\qquad k=0,\ldots,r=\dim\leaf
\]
And we can write all these equations in matrix form
\[
\left(\begin{array}{cccc}1 & 1 & \cdots & 1 \\\lambda_1 & \lambda_2 & \ldots & \lambda_r \\\vdots & \vdots & \ddots & \vdots \\\lambda_1^r & \lambda_2^r & \cdots & \lambda_r^r\end{array}\right)\cdot\left(\begin{array}{c}a_1^2 \\a_2^2 \\\vdots \\a_r^2\end{array}\right)=\left(\begin{array}{c}c_1 \\c_2 \\\vdots \\c_r\end{array}\right)
\]
The matrix on the left is a Vandermonde matrix, and since the eigenvalues are distinct, is invertible. In particular the functions $a_i^2$ are all constant.
\qed

%------------------------- 2 DIM SRF ON SPHERES --------------------------------------- 2 DIM SRF ON SPHERES --------------------------------

\bigskip

\section{2-dimensional SRF on spheres}\label{2dim}

\bigskip

In this section we prove that 2-dimensional foliations are homogeneous. Let $(S^n,\fol)$ be a 2-dimensional SRF, which we assume to be without 0-dimensional leaves by Proposition \ref{split}. By Proposition \ref{sigma1} $\Sigma_1$ is a union of $k$ disjoint totally geodesic spheres. Let $S$ be one such component, and $S'$ be the sphere at distance $\pi/2$ from $S$. Both $S$ and $S'$ are saturated by leaves, $S^n=S\star S'$, and $\fol|_{S'}$ contains the remaining $k-1$ components of $\Sigma_1$.

Notice that $k>0$ since otherwise $\fol$ would be regular and 2 dimensional, which is not possible in a sphere. For the same reason, $k>1$ otherwise $\fol|_{S'}$ would be regular and 2 dimensional as well.  Also, if $k=2$ then $\fol|_{S'}$ is 1-dimensional. In fact, if it weren't so, then by replacing $S'$ with $S'$ in the discussion above we could write $S'=S''\star S'''$, as a join of spheres where $\fol|_{S''}$ is 1-dimensional, and $\fol|_{S'''}$ is regular and 2-dimensional.

We will now prove the result by induction on $k$. More precisely, we will prove the following:
\begin{itemize}
\item If $k=2$, then $(S^n,\fol)=(S,\fol|_S)\star(S',\fol|_{S'})$. Since $\fol|_S$, $\fol|_{S'}$ are 1-dimensional, they are both homogeneous and hence so is  $\fol$. Moreover, the group acting on $S^n$ is $\RR^2$.
\item If $k>2$ and  $\fol|_{S'}$ is homogeneous by the action of $\RR^2$, so is $\fol$.
\end{itemize}

\subsection*{Proof of (1).}

Any horizontal geodesic from $S$ to $S'$ passes through 2-dimensional leaves, and meets 1 dimensional leaves at the endpoints. The rigidity condition as in point (3) of Proposition \ref{usefulprop} are met, and therefore we have the splitting $\fol=\fol|_{S}\star \fol|_{S'}$. Since $\fol|_S$, $\fol|_{S'}$ are  one dimensional, they are given by orbits of $\RR$-actions $\rho_S,\,\rho_{S'}$, and therefore $\fol$ is given by orbits of the $\RR^2$-action $\rho_S\times\rho_{S'}$.

\subsection*{Proof of (2).}

Let $\rho: \RR \curvearrowright S$, $\rho':\RR^2\curvearrowright S'$ be the actions on $S$ and $S'$, respectively. Notice first that since we are assuming that $k>2$, $S'$ contains 2-dimensional leaves. Take $p\in S$, $p'\in S'$ a regular point, and consider the unit length (horizontal) geodesic $\gamma$ from $p'$ to $p$. Given $v\in \RR^2$ set $v^*\in \V_{p'}$ the action field corresponding to $v$, at the point $p'$. Since $p'$ is at a regular point, there exists a unique holonomy Jacobi field $J_v$ along $\gamma$ such that $J_v(0)=v^*$. This Jacobi field is always vertical, and in particular $J_v(\pi/2)\in \V_p$ is given by some action field $w^*$ corresponding to a vector $w\in \RR$. Consider the map $\pi:\RR^2\to \RR$ given by $v\lmt w$ obtained as we just described: this is a linear map, and the (isometric) linear action $\rho'\oplus \rho\circ \pi:\RR^2\curvearrowright S^n$ gives rise to a homogeneous SRF $\fol'$. By definition, $\fol|_S=\fol'|_S$, $\fol|_{S'}=\fol'|_{S'}$, and $\fol|_{\gamma}=\fol'|_{\gamma}$. We want to show that $\fol=\fol'$, which will prove that $\fol$ is homogeneous. The way we want to prove this is by showing that both foliations $\fol, \fol'$ are determined by their restrictions on $S$, $S'$ and $\gamma$. It is enough to show that $\fol$ is uniquely determined by its restrictions as just said.

First of all, $\fol$ is uniquely determined on $C_{\theta}=\cos(\theta)\leaf_{p'}+\sin(\theta)\leaf_p$, for any fixed $\theta$. In fact, $C_{\theta}$ is isometric to $\leaf_{p'}\times \leaf_{p}$ and hence is flat (both $\leaf_p$ and $\leaf_{p'}$ are flat). Furthermore, $C_{\theta}$ is 3-dimensional, and it is saturated by 2-dimensional flat leaves. Lifting $\fol$ to $\wt{\fol}$ in the universal cover $\wt{C}_{\theta}=\RR^3$, $\wt{\fol}$ is a foliation by parallel planes and therefore it is uniquely determined by the vertical space at a single point. Then the same is true for $\fol|_{C_{\theta}}$, which is determined by the vertical space $\V_{\cos(\theta)p'+\sin(\theta)p}=\V_{\gamma(\theta)}$.

Now consider any other point $\ol{q}=\cos(\sigma)q'+\sin(\sigma)q$, where $q'\in S'$ is again regular, and $q\in S$. Again, set $C'_{\sigma}=\cos(\sigma)\leaf_{q'}+\sin(\sigma)\leaf_q$, which is 3-dimensional and is saturated by 2-dimensional leaves.

Let $c$ be the minimizing geodesic between $\ol{q}$ and $C_{\theta}$, and $\ol{p}$ be the end point of $c$. Then there is a neighborhood $U$ of $\ol{p}$ in $\leaf_{\ol{p}}$ such that $f(x)=dist(U, x)$ is regular around $\ol{q}$, and the level sets of $f$ give rise to a codimension 1 foliation in $C'_{\sigma}$. By definition of SRF this foliation coincides with $\fol$ around $\ol{q}$, and therefore $\fol$ is determined on an open dense set of $S^n$. In particular, $\fol=\fol'$ on an open dense set of $S^n$, and by Lemma 4.4 of \cite{LT} they coincide everywhere.
\qed

%-------------------------- 3DIM SRF ON SPHERES ------------------------------------ 3 DIM SRF ON SPHERES ---------------------

\bigskip

\section{3-dimensional SRF on spheres: subcases}\label{cases}

\bigskip

Let $(S^n,\fol)$ a SRF with $3$-dimensional regular leaves and no 0-dimensional leaf. If $p$ is a point of a regular leaf $\leaf$, we can subdivide $\H_p$ into
\[
\H_p=\H_{0,p}\cup\H_{1,p}\cup\H_{2,p}\cup\H_{3,p}\qquad\textrm{where}\quad \H_{i,p}:=\left\{x\in \H_p\st rank(A_x)=i\right\}
\]
Notice that if $X$ is a basic vector field on some open set $U$, and $p,q\in U$ belong to the same leaf, then $X_p\in \H_{i,p}$ iff $X_q\in \H_{i,q}$ . Therefore, the decomposition of $\H_p$ only depends on the leaf, and not on the specific point.
\\

The discussion about 3-dimensional foliations will be divided into 5 cases:
\begin{itemize}
\item[I.] There exists a point $p$ and some $x\in\H_{0,p}$, such that $S_x$ has 3 different eigenvalues.
\item[II.] There exists a point $p$ and some $x\in\H_{0,p}$, such that $S_x$ has 2 different eigenvalues.
\item[III.] For some point $p$, $\H_{0,p}=\emptyset$, $\H_{1,p},\H_{2,p}\neq\emptyset$.
\item[IV.] For every regular point $p$, $\H_{0,p}=\H_{1,p}=\emptyset$.
\item[V.] For every regular point $p$, $\H_{0,p}=\H_{2,p}=\emptyset$.
\end{itemize}

\begin{prop}
There are no other cases to study.
\end{prop}
\proof First of all, notice that if there is any horizontal vector $x$ such that $A_x$ has rank 3, then this foliation is \emph{substantial} and therefore homogeneous by the work of Grove and Gromoll \cite{GG}. We can then suppose that for every horizontal vector $x$, $A_x$ has rank at most 2. Cases I and II cover all the possible cases in which sone $x$ ha rank 0, since we cannot have an $x$ such that $A_x=0$ and $S_x$ has 1 eigenvalue $\lambda$. In fact, if this happened, then by lemma \ref{usefullem} the geodesic $\exp(tx)$ would meet a 0-dimensional singular leaf at time $t=1/tan(\lambda)$, which contradicts the assumption of having no 0-dimensional leaves.

Suppose now that every $A_x$ has rank 1 or 2, and take a regular point $p$. If we are not in case III, then either every $A_x$, $x\in \H_p$ has rank 1, or every $A_x$, $x\in \H_p$ has rank 2. Suppose they all have rank 1, and take a horizontal geodesic $\gamma$ from $p$ to any other regular point $q$. By the formula $$(\nabla^v_xA)_xy=2S_xA_xy$$ (see \cite{GG}, page 149) the rank of $A_{\gamma'(t)}$ is constant, therefore $\rk(A_{\gamma'(1)})=\rk(A_{\gamma'(0)})=1$. In particular there exists a vector $\gamma'(1)\in \H_q$ whose corresponding O'Neill tensor has rank 1.  We then have two possibilities:
\begin{itemize}
\item There exists another point $q$ that falls in case III.
\item For every point in the regular stratum $\H=\H_1$, and this corresponds to case V.
\end{itemize}
\qed
% ------------------------------------------------------- CASE 1 ----------------------------------------- CASE 1 --------------------------------

\bigskip

\section{3-dimensional SRF on spheres: case I}

\bigskip

Suppose there is a regular leaf $\leaf_0$ along which we have a basic vector field $X$ which is parallel (meaning $A_X^*=0$), and such that $S_X$ has 3 different eigenvalues $\lambda_i$, $i=1,2,3$. Let $E_i$, $i=1,2,3$ be the eigendistributions. Because of the Codazzi equation, these distributions are autoparallel, and hence integrable and totally geodesic. In fact, let $V_1,V_2\in E_i$ and $W\in E_j$ with $i\neq j$. Then by the Codazzi equation:
\[
\scal{\left(\nabla_{W}S\right)_XV_1,V_2}=\scal{\left(\nabla_{V_1}S\right)_XW,V_2}
\]
The left hand side is equal to
\[
\scal{\nabla_{W}(S_XV_1),V_2}-\scal{S_X\nabla_{W}V_1,V_2}=\lambda_i\scal{\nabla_WV_1,V_2}-\lambda_i\scal{\nabla_WV_1,V_2}=0
\]
The right hand side is equal to
\begin{eqnarray*}
\scal{\nabla_{V_1}(S_XW),V_2}-\scal{S_X\nabla_{V_1}W,V_2}&=&\lambda_j\scal{\nabla_{V_1}W,V_2}-\lambda_i\scal{\nabla_{V_1}W,V_2}\\&=&(\lambda_i-\lambda_j)\scal{W,\nabla_{V_1}V_2}
\end{eqnarray*}
Therefore $(\lambda_i-\lambda_j)\scal{W,\nabla_{V_1}V_2}=0$, and since $\lambda_i\neq \lambda_j$ we obtain $\nabla_{V_1}V_2\perp E_j$. Since $E_j$ was chosen arbitrarily, $\nabla_{V_1}V_2\in E_i$.

Since our eigendistributions are 1-dimensional, we can locally find vertical vector fields $V_i$, $i=1,2,3$, such that:
\begin{itemize}
\item $\scal{V_i,V_j}=\delta_{ij}$.
\item $\nablal_{V_i}V_i=0$.
\item $S_XV_i=\lambda_iV_i$.
\end{itemize}
Where $\nablal$ denotes the Levi-Civita connection on $\leaf$.

Call $\E$ the vector space generated by the vector fields $V_i$.

We will now prove that $\E$ satisfies the conditions of Theorem \ref{HomThm}, thus proving that the foliation is homogeneous. The proof will proceed through the following steps:
\begin{enumerate}
\item The $V_i$ defined above are eigenvalues for any shape operator of basic fields, and in particular $S_X\mc{E}=\mc{E}$ for any basic $X$.
\item $\mc{E}$ is preserved under the usual Lie Brackets, and in particular is a (finite dimensional) Lie algebra.
\item $A_XY\in \mc{E}$ for any basic $X,Y$.
\end{enumerate}
\subsection*{Proof of (1).}
We use the Ricci equation $$\scal{\nabla^h_U\nabla^h_VX-\nabla^h_V\nabla_U^hX-\nabla^h_{[U,V]}X,Y}=\scal{[S_X,S_Y]U,V}$$
where again $A^*_X=0$. All the three terms on the left hand side are zero: in fact, if $V$ is any vertical vector field, $\nabla^h_VX=\nabla_X^hV=-A_X^*V=0$.
Therefore $[S_X,S_Y]=0$ for every $Y$ and since $S_X$ has 3 different eigenvalues then the $V_i$ defined above are eigenvalues for any shape operator. In particular we have functions $\lambda_i:\H_p\ra\RR$ such that
\[
S_YV_i=\lambda_i(Y)V_i.
\]
Of course $\lambda_i=\lambda_i(X)$. This tells us that the vectors $B(V_i,V_i)$ are basic, that $B(V_i,V_j)=0$, and that the curvatures $k_{ij}=sec(V_i,V_j)$ are constant.

\subsection*{Proof of (2).}
Notice that
\begin{eqnarray*}
\scal{\nablal_{V_i}V_j,V_i}&=&-\scal{V_j,\nablal_{V_i}V_i}=0\\
\scal{\nablal_{V_i}V_j,V_j}&=&V_i\left(-{1\over2}\|V_j\|^2\right)=0
\end{eqnarray*}
so $\nablal_{V_i}V_j$ is either $0$ if $i=j$, or is a multiple of $V_k$, where $k$ is the third possible index
\[
\nablal_{V_i}V_j=f_{ij}V_k
\]
Deriving the equation $\scal{V_i,V_j}=\delta_{ij}$ we obtain
\begin{equation}\label{skew}
0=V_k\scal{V_i,V_j}=\scal{\nablal_{V_k}V_i, V_j}+\scal{V_i,\nablal_{V_k}V_j}=f_{ki}+f_{kj}
\end{equation}
that means $f_{12}=-f_{13}$, $f_{23}=-f_{21}$, $f_{31}=-f_{32}$. 
Now consider the equation
\begin{eqnarray*}
f_{ij}f_{ji}&=&\scal{\nablal_{V_i}V_j,\nablal_{V_j}V_i}=\\
&=&-\scal{\nablal_{V_j}\nablal_{V_i}V_j,V_i}\\
&=&-\scal{R^{\leaf}(V_j,V_i)V_j,V_i}-\scal{\nablal_{V_i}\nablal_{V_j}V_j,V_i}-\scal{\nablal_{[V_j,V_i]}V_j,V_i}\\
&=&k_{ij}-(f_{ji}-f_{ij})\scal{\nablal_{V_k}V_j,V_i}\\
&=&k_{ij}+(f_{ij}-f_{ji})f_{kj}
\end{eqnarray*}
Together with the equivalences among the functions $f_{ij}$ found above, we get
\begin{eqnarray*}
-f_{12}f_{23}=k_{12}-f_{31}(f_{12}+f_{23})\\
-f_{23}f_{31}=k_{23}-f_{12}(f_{23}+f_{31})\\
-f_{31}f_{12}=k_{31}-f_{23}(f_{12}+f_{23})
\end{eqnarray*}
If we denote $b_1=f_{31}f_{12}$, $b_2=f_{12}f_{23}$, $b_3=f_{23}f_{31}$, then the system above becomes a system of linear equations in the indeterminates $b_i$, which has a unique (constant!) solution. 
From this, it follows immediately that the $f_{ij}^2={b_ib_j\over b_k}$ are constant, and $\span\{V_1,\,V_2,\,V_3\}$ is a 3-dimensional Lie algebra.

\subsection*{Proof of (3).}
We know that $S_X$ has 3 different eigenvalues for the special $X$, and in particular there are 3 different eigenvalues for almost all basic vector fields on $\leaf$. Also, the eigenvectors are $V_1,V_2,V_3$ defined before. Then by Proposition \ref{Slem} the inner product $\scal{A_ZY,V_i}$, $i=1,2,3$, is constant for all basic $Y$ and all basic $Z$ with 3 eigenvalues. In particular, $A_Z^*V_i\in\B$, and therefore $A_Z\E\In \B$. By continuity, $A_Z\E\In \B$ for all basic $Z$.
\qed

% ------------------------------------------------------- CASE 2 ----------------------------------------- CASE 2 --------------------------------

\bigskip

\section{3-dimensional SRF on spheres: case II}

\bigskip

In this section we suppose we are in case II, i.e. there is a leaf $\leaf_0$, with a horizontal vector $x\in \H_p$, $p\in \leaf_0$ such that $A_x^*=0$ and $S_x$ has 2 different eigenvalues.

Again we want to prove that $\fol$ is homogeneous. The strategy will be the following:
\begin{enumerate}
\item We prove that there exists a global Killing vector field that preserves the foliation, whose flow has a nonempty fixed point set $N$, and $N$ is a totally geodesic sphere $S^k\In S^n$.
\item We show that $N=S^k$ is preserved by the foliation, and therefore the orthogonal sphere $S^{n-k-1}$ is preserved as well.
\item If $\fol|_N$ has 1-dimensional leaves, we show that we can reduce orselves to case I.
\item If $\fol|_N$ has no 1-dimensional leaves, we prove that $N=S^2$ with $\fol|_{N}$ consisting of only one leaf, and the foliation splits as $(S^{n-k-1},\fol_1)\star (S_2,\fol|_{N})$
\end{enumerate}

\subsection*{Proof of (1).}
Let $X$ be the basic vector field such that $X_p=x$, let $\lambda_1$, $\lambda_2$ the different eigenvalues of $S_X$, and let $E_1, E_2$ be the corresponding eigen-distributions, such that $\dim E_i=i$. Locally around a point $q$, let $W$ be an unit length vector field spanning $E_1$, and $V_1, V_2$ an orthonormal frame of $E_2$.  We already know (cfr. for example \cite{CO}, remark 1.3) that the eigen-distributions are parallel and the following equations hold:
\begin{equation}\label{Wkilling}
\nablal_WW=0,\quad \nablal_W V_i=0.
\end{equation}
Notice that along $\gamma(t):=\exp_ptX$, the Jacobi field with initial value $W$ will have a zero when meeting a 2 dimensional leaf $\leaf_2$ at distance $d_1={1/\tan\lambda_1}$. Similarly the holonomy Jacobi fields with initial values in $E_2$ will have a zero when meeting a one dimensional leaf $\leaf_1$ at distance $d_2=1/\tan\lambda_2$.

From the Ricci equation, $W$ is an eigenvalue of every shape operator. Moreover, since $W$ is the unique eigenvector of its eigenvalue, by lemma \ref{Slem} $\scal{A_XY,W}$ is constant for all basic $X,Y$, and in particular $A_X^*W$ is basic for all basic $X$. As in the proof of the Homogeneity Theorem \ref{HomThm}, we argue that $W$ is the restriction to $\leaf_0$ of a global Killing field on $S^n$. In fact, let $\phi^t$ be the flow of $W$ for time $t$: because of the equations \ref{Wkilling}, we know that $W$ is a Killing field along $\leaf_0$, and $\phi^t$ is an isometry. Moreover $W$ commutes with $V_1,V_2$ and thus $\phi^t_*V_i=V_i\circ \phi^t$. Finally, define $\Phi^t:\nu(\leaf_0)\ra \nu(\leaf_0)$ by $$\Phi^t(X_p)=X_{\phi^t(p)}$$ for every basic vector field $X$.

We now apply the Fundamental Theorem of Submanifold Geometry:

\begin{teor}[Fundamental Theorem of Submanifold Geometry]
Let $M,N\In S^n$ be two submanifolds, let $f:M\ra N$ be an isometry and $F:\nu(M)\ra \nu(N)$ a bundle map over $f$. Suppose that the couple $(F,f)$ preserves the shape operator and the normal connection, i.e. for every sections $X\in \Gamma\big(\nu(M)\big)$ and $V\in \vf(M)$ the following equations are satisfied: $$f_*(S^M_XV)=S^n_{F(X)}(f_*V)\qquad F\big(\nabla^{\perp}_VX\big)=\nabla^{\perp}_{f_*V}\big(F(X)\big).$$
Then there exists a global isometry $h:S^n\ra S^n$ such that $h|_M=f$, $h_*|_{\nu(M)}=F$.
\end{teor}
In our case, $M=N$ is an orbit of $W$, and $TM$ is spanned by $W$, while $\nu(M)=\nu(\leaf)\oplus \span\{V_1,V_2\}$. Furthermore $f=\phi^t|_{M}$ and $F$ is given by $$F=\left(\phi^t_*\right)|_{\span\{V_1,V_2\}}\oplus \Phi^t.$$ We now check that the hypotheses of the fundamental theorem are satisfied:
If $X$ is basic, $S^M_XW=\scal{S_XW,W}W=S_XW=\lambda_1(X) W$ for some constant $\lambda_1(X)$, and therefore
\begin{eqnarray*}
f_*(S^M_XW)&=&\lambda_1(X)\phi^t_*W=\lambda_1(X) W\circ \phi^t=\\
&=&\left(S_XW\right)\circ \phi^t=S_{X\circ \phi^t}(W\circ\phi^t)=S_{F(X)}(f_*W)
\end{eqnarray*}
This proves the first condition, for $X$ basic.

Moreover, $$S^M_{V_i}W=\scal{\nablal_WV_i,W}W=0$$ and this proves the first condition, for the vectors $V_i$.

If $X$ is basic, then $\nabla^{\H}_{W}X=-A_X^*W$ is basic, and  $$ \nabla_W^{\perp}X=\nabla^{\H}_WX+\sum_{i=1}^2\scal{\nabla_WX,V_i}V_i=\nabla^{\H}_WX+\sum_{i=1}^2\scal{S_XW,V_i}V_i=\nabla^{\H}_WX$$ Therefore
\begin{eqnarray*}
\Phi^t\big(\nabla^{\perp}_WX\big)=\Phi^t\big(\nabla^{\H}_WX\big)&=&-\left(A_X^*W\right)\circ\phi^t\\ &=&-A_{X\circ\phi^t}^*(W\circ\phi^t)\\&=&-A^*_{\Phi^tX}(\phi^t_*W)=\nabla^{\perp}_{\phi^t_*W}\Phi^tX.
\end{eqnarray*}
This proves the second condition, for $X$ basic. Moreover $$\nabla^{\perp}_{W}V_i=\sum_{i=1}^2\scal{\nablal_{W}V_i,V_i}V_i+\nabla^{\H}_WV_i=0$$ and this proves the second condition, for the vectors $V_i$.

Therefore, $\phi^t$ extends to a global isometry $\psi^t$, and since $\psi^t_*|_{\nu(M)}=\Phi^t$ preserves basic fields, then $\psi^t$ preserves the foliation. Thus the Killing vector field $K={d\over dt}|_{t=0}\psi_t$ restricts to $W$ on $\leaf_0$, and to a holonomy Jacobi field along horizontal geodesics. As a holonomy Jacobi field, it cannot be zero on the regular leaves, so the fixed point set $N:=Fix(\RR)$ only intersects singular leaves. 

\subsection*{Proof of (2).}
Notice first of all, since the $\RR$ action preserves the leaves, it acts isometrically on each of them. In particular, a one dimensional leaf is either acted on transitively, or it is totally fixed, and this proves that $\fol|_N$ preserves one dimensional leaves.

For two dimensional leaves, notice that if the leaf is not fixed then the fixed point set consists of a discrete set of points.
Now, we know that the singular two dimensional leaf $\leaf_2$ described above is totally contained in the fixed point set. Suppose that another leaf intersects $N$ at some point $q$. Then $\RR$ acts almost effectively on $T_{q}\leaf_q$ by isometries. Let $c:[0,1]\ra N$ a horizontal geodesic from $q$ to $q'\in \leaf_1$, and let $J$ a holonomy Jacobi field along $c$, such that $J(0)=v\in T_q\leaf_q$, $J(1)=w\in T_{q'}\leaf_1$. If $t\in\RR$ is an element acting on $S^n$, on one hand it preserves $c$, and sends $J$ to $t_*J$ still a Jacobi field. On the other hand, since it preserves the foliation, we know that $t_*J$ is also a holonomy Jacobi field. But $t_*J(0)=t_*v\neq v$, while $t_*J(1)=t_*w=w$. It follows that $(J-t_*J)$ is a nonzero holonomy Jacobi field that goes to zero at some ``regular leaf'' (in the singular stratum), and this can't be true. In particular, if a 2-dimensional leaf intersects $N$, it must be contained in it.
The result follows since we observed that no 3-dimensional leaves intersect $N$.
\\

$N$ is a totally geodesic sphere, and by the lemma above we know that it has a 2-dimensional SRF, hence it is either $N=S^2$ with only one leaf, or it is a homogeneous foliation under the action of $\RR^2$.

\subsection*{Proof of (3).}
If $N\neq S^2$, then it contains at least 2 connected components of $\Sigma_1$, which we denote $C_2$ and $C_3$. Moreover, $\leaf_1$ is contained in a third component, call it $C_1$. The horizontal geodesic $\gamma(t)=\exp tx$ defined above meets $\leaf_1$ in $p_1$ and $\leaf_2$ in $q_1$. Without loss of generality, we can assume that $q_1$ can be written as $\cos(\theta)p_2+\sin(\theta)p_3$, for some $\theta\in (0,\pi/2)$, where $p_2\in C_2,\, p_3\in C_3$. Take now the 2-sphere spanned by $p_1,p_2,p_3$, and consider the geodesic $c$ passing through $q_1$ and $q_2=\cos(\alpha)p_1+\in(\alpha)p_2$, for some $\alpha\in (0,\pi/2)$. The geodesic $c$ satisfies the following properties:
\begin{itemize}
\item It is horizontal, because it is horizontal at $q_1$. In fact $c'(0)$ is given by a combination of $\gamma'(d_1)$ and $s'(\theta)$, where $s(t)=\cos(t)p_2+\sin(t)p_3$ is horizontal.
\item For an appropriate chice of $\alpha$, $c$ passes through 3-dimensional leaves.
\item It meets $\Sigma_2$ at 6 points, namely $q_1, q_2$, some $q_3=\cos(\beta)p_1+\sin(\beta)p_3$, and their antipodal points.
\end{itemize}
Therefore, taking a regular point $c(\epsilon)$, the vector $c'(\epsilon)$ satisfies the hypotheses of case I, and the foliation is homogeneous.

\subsection*{Proof of (4).}
If $N=S^2$, then regular leaves have (common) universal cover $S^2\times \RR$.  Let $N'=S^{n-3}$ be the totally geodesic sphere at distance $\pi/2$ from $N$. Since $N$ is saturated by leaves, so is any distance tube around $N$, and in particular $N'$ is saturated as well. $N'$ is acted on by the isometries $\psi^t$ without fixed points, and this gives a regular 1-dimensional foliation on $N'$, whose leaves are contained in our original SRF. We are going to show that the foliation on $N'$ actually consists of this 1-dimensional foliation. To do so, first notice that $N'$ does not contain 2-dimensional leaves. In fact, if there were one such a leaf $\leaf_2$, there would be a submersion $\leaf_{reg}\to \leaf_2$ that would induce a fibration $F^1\to S^2\times \RR\to \leaf_2$ with one dimensional fiber. Looking at the exact sequence in homotopy, we obtain $\pi_2(\leaf_2)=\ZZ$ which implies $\leaf_2$ is diffeomorphic to either $S^2$ or $\RR\PP^2$. But in either case, $\leaf_2$ cannot admit a regular 1-dimensional foliation. In particular, there is only one component of $\Sigma_1$ in $N'$. In fact, if $C_1$, $C_2$ were two components of $\Sigma_1$, then by Proposition \ref{sigma1} the join $C_1\star C_2$ would contain 2-dimensional leaves. Call $N_1$ the component of $\Sigma_1$ in $N'$. Then the sphere $N_2\In N'$ at distance $\pi/2$ from $N_1$ preserves the foliation as well, and does not contain any 1-dimensional leaf. It has then to contain only 3-dimensional leaves. The foliation in $N_2$ is then regular, and by the result of Grove and Gromoll in \cite{GG} it is homogeneous. In particular the leaves are covered by spheres, but this contradicts the fact that regular leaves are covered by $S^2\times \RR$. As a consequence $N_2=\emptyset$ and $N'=N_1$ only contains 1-dimensional leaves, which is what we wanted to prove.
By Proposition \ref{usefulprop}, this implies that the foliation on $S^n$ is a join $$(S^n,\fol)=(N,N)\star (N',\psi^t).$$
\qed

% ------------------------------------------------------- CASE 3 ----------------------------------------- CASE 3 --------------------------------

\bigskip

\section{3-dimensional SRF on spheres: case III}

\bigskip

In this section we assume we are in case III, i.e. there is a regular point $p\in S^n$ such that for any $x\in \H_p$, $A_x$ has either rank 1 or 2. The strategy in this case consists of 2 main steps:
\begin{enumerate}
\item First we show that $\A=\left\{A_xy\st x,y\in \H_p\right\}$ is 2-dimensional, and therefore the orthogonal space $\V_p\cap \A$ is locally spanned by one vector field $U$.
\item Finally, we show that the sheaf on $\leaf_p$ locally consisting of $U$ and $\span \{A_XY\}$, with $X,Y$ basic, satisfies the condition of the homogeneity Theorem \ref{HomThm}. Here we heavily rely on the fact that, $\A$  2-dimensional by the previous point, then for almost every $x\in \H_p$, $\im(A_x)=\A$.
\end{enumerate}

\subsection{Step (1).} Remember we defined $\H_1:=\left\{x\in \H_p\st \rk(A_x)=1\right\}$. Call $\X=\PP(\H_1)\In\PP^{h-1}$, and decompose $\X$ into its irreducible components
\[
\X=\X_1\cup\ldots\cup \X_p.
\]
Also, define Zariski open subsets $\Y_1,\ldots \Y_p\in\H_2$
\[
\Y_i=\left\{y\in\H_2\st A_y\X_i\neq 0\right\}
\]
\begin{lem}
Define the sets in $\V$:
\[
\Z_i=\cup_{x\in \X_i}\im\{A_x\}
\]
Then the following are true:
\begin{itemize}
\item $\Y_i\neq0$ for all $i$, and hence they are all open dense sets in $\H_2$.
\item For every $y\in \Y_i$, $\Z_i\In \im A_{y}$.
\end{itemize}
\end{lem}
\proof\begin{itemize}
\item Suppose $\Y_i=\emptyset$. Then for all $y\in \H_2$, and all $x\in \X_i\In \H_1$, $A_yx=0$. But since $\H_2$ is open and dense in $\H$, we obtain that $A_{\H}x=A_x\H=0$, and this means $x\in \H_0$, contradicting $x\in\H_1$.
\item Take $y\in \Y_i$, and $x\in \X_i$ such that $A_yx\neq0$. Since $\X_i$ is irreducible, $\X_i-\ker A_y$ is open and dense in $\X_i$. Now, the diagram below commutes,
\[\begindc{\commdiag}[50]
\obj(0,0)[PZ]{$\PP(\X_i-\ker A_y)$} \obj(0,1)[PZ']{$\PP(\X_i)$} \obj(3,1)[PV]{$\PP(\V)$} \mor{PZ}{PZ'}{}[1,5] \mor{PZ}{PV}{$A_y$} \mor{PZ'}{PV}{$[x]\lmt[\im\{A_x\}]$}
\enddc\]
And in particular every $x\in \X_i-\ker A_y$ gets mapped to $\im\{A_y\}$. By continuity, the whole $\X_i$ gets mapped in the same plane.
\end{itemize}
\qed
As a corollary, we have the following
\begin{lem}
There is a plane $\pi$ such that $\im\{A_y\}=\pi$ for all $y\in \H_2$.
\end{lem}
\proof There are two possibilities: either the $\X_i$ get mapped to the same line, or not.
Suppose first that not all $\X_i$ get mapped to the same line. Then for every $y\in \Y_1\cap\ldots\cap \Y_p$ the plane $\im\{A_y\}$ contains all the $\Z_i$. But then all these planes must be the same plane $\pi$, and since $\Y_1\cap\ldots\cap \Y_p$ is dense in $\H$, then the image of every $A_x$ lands in $\pi$.
\\

Suppose now that every $\X_i$ gets mapped to the same line $L$, and for every $y\in \H_2$, the plane $\im\{A_y\}$ contains $L$. Suppose there are $y,y'\in \H_2$ such that $\im\{A_y\}\neq \im\{A_{y'}\}$. The set
\[
A_{y}^{-1}(\V-L)\cap A_{y'}^{-1}(\V-L)\In \H_2
\]
is nonempty (and in particular is open and dense). If $x$ belongs to that set, then $A_xy$, $A_xy'$ span a plane that does not contain $L$, and this is a contradiction. in particular, $\im\{A_y\}=\im\{A_{y'}\}$.
\qed
In particular, there is a unit length vector field $U$ such that $A_X^*U=0$ for every $X$. At the same time, $U^{\perp}$ is spanned by $A_XY$, with basic $X,Y$.
\subsection{Step (2).} Define
\begin{eqnarray*}
\B&=&\{\textrm{basic vector fields}\}\\
\A&=&\left\{A_XY\st X,Y\in\B\right\}\\
\C&=&\scal{U}
\end{eqnarray*}
and $\mc{E}:=\A\oplus\C$. The goal of this part is to prove the following proposition
\begin{prop}\label{caseIII}
The following properties hold:
\begin{enumerate}
\item $\scal{V_1,V_2}=cost.$ for all $V_1,V_2\in\mc{E}$.
\item $S_X\mc{E}\In\mc{E}$ for every basic $X$.
\item $[\A,\A]\In \mc{E}$.
\item $[\A,\C]\In \mc{E}$.
\end{enumerate}
In particular $\mc{E}$ is a Lie algebra.
\end{prop}

\subsection*{Proof of (1).}

Say $V_1=A_{X_1}X_2+aU$, $V_2=A_{X_3}X_4+bU$, where $a,b$ are constants, and $X_1,\ldots, X_4$ are basic. Then
\[
\scal{V_1,V_2}=\scal{A_{X_1}X_2,A_{X_3}X_4}+ab
\]
So it remains to show that $\scal{A_{X_1}X_2,A_{X_3}X_4}$ is constant.

It is enough to prove it for a dense set of possible $X_1,\ldots X_4$, so we can suppose that $X_1\in \H_2$, and that $A_{X_1}X_3, A_{X_1}X_4$ span $\im\{A_{X_1}\}$. One can replace $X_3$, $X_4$ by linear combinations $\bar{X}_3, \bar{X}_4$ such that $A_{X_1}\3,A_{X_1}\4$ is an orthonormal basis.
In particular
\begin{eqnarray*}
A_{X_1}X_2&=&\scal{A_{X_1}X_2,A_{X_1}\3}A_{X_1}\3+\scal{A_{X_1}X_2,A_{X_1}\4}A_{X_1}\4\\
&=&c_1A_{X_1}\3+c_2A_{X_1}\4
\end{eqnarray*}
and
\begin{eqnarray*}
\scal{A_{X_1}X_2,A_{\3}\4}=c_1\scal{A_{X_1}\3,A_{\3}\4}+c_2\scal{A_{X_1}\4,A_{\3}\4}
\end{eqnarray*}
is constant.

\subsection*{Proof of (2).}

First of all, we show that
\[
S_X\A\In\mc{E}
\]
Again it is enough to prove it for a dense subset of $X$, so we can suppose $X\in \H_2$, in particular, by what we proved before, every $V\in \A$ can be written as $A_XY$ for some $Y$. Pick an orthonormal basis in $\mc{E}$, $\{A_XY_1,A_XY_2,U\}$, and set
\[
S_XA_XY=f_1A_XY_1+f_2A_XY_2+f_3U.
\]
We know immediately that
\begin{eqnarray*}
f_1&=&\scal{S_XA_XY,A_XY_1}\\
f_2&=&\scal{S_XA_XY,A_XY_2}
\end{eqnarray*}
are constant, by Lemma \ref{Slem}. Also, by using Lemma \ref{Slem} again we know that
\[
f_1^2+f_2^2+f_3^2=\scal{S_XA_XY,S_XA_XY}=\scal{S^2_XA_XY,A_XY}
\]
is constant as well, hence $f_3$ is constant too. In particular, $S_XA_XY\in\mc{E}$.

We now prove that $S_XU\in \mc{E}$, in other words that
\[
\scal{S_XU,A_XY},\;\scal{S_XU,U}
\]
are constant. We already know that
\[
\scal{S_XU,A_XY}=\scal{S_XA_XY,U}
\]
is constant, so we only need to show that $\scal{S_XU,U}$ is constant. Suppose first that $S_X$ has 3 eigenvalues, with eigenvectors $V_1,V_2,V_3$. Take the othonormal basis $A_XY_1,A_XY_2,U$ as before. We know by Proposition \ref{Slem} we know that $A_XY_1$ $A_XY_2$ can be written as combinations of the $V_i$'s with constant coefficients, hence so it must be $U$,
\[
U=\sum_{i=1}^3 a_iV_i.
\]
Then
\[
\scal{S_XU,U}=\sum_{i=1}^3a_i^2\lambda_i
\]
is constant.

\subsection*{Proof of (3).}

This part is taken directly from \cite{GG}. First of all, we show that $[\A,\A]\In[\A,\B]\oplus[[\A,\B],\B]$: if $X,Y\in\B$, and $T\in\A$,
\begin{eqnarray*}
2[A_XY,T]&=&[[X,Y]^v,T]=[[X,Y],T]-[[X,Y]^h,T]\\
&=&-[[Y,T],X]-[[T,X],Y]-[[X,Y]^{h},T]
\end{eqnarray*}
Now we are going to prove that $[\A,\B]\In\mc{E}$: for this, remember once more that if $X\in\H_2$, every $A_YZ$ can be written as $A_XY'$ for some $Y'$. Hence
\begin{eqnarray*}
[X,A_YZ]=[X,A_XY']&=&\nabla^v_XA_XY'+S_XA_XY'\\
&=&3S_XA_XY'+A_X\nabla_X^hY'-A_{Y'}\nabla^h_XX
\end{eqnarray*}
that implies $[\A,X]\In \A\oplus S_X\A\In \mc{E}$.

The last thing to prove now is that $[\B,[\B,\A]]\In \mc{E}$, and for this we only have to show that $[S_XA_XY,Z]\in \mc{E}$, for basic $X,Y,Z$. Now, this would be obvious if $S_X\A\In \A$ for all $X\in \B$, so we can suppose that for almost every $X$, $S_X\A\neq \A$. Here we can work with $X,Y,Z$ generic basic vector fields, so i suppose that $S_X\A\neq \A$, $S_Z\A\neq \A$. But since $\A\In \mc{E}$ is a plane in 3-space, we can write
\[
S_XA_XY=S_ZA_ZY_1+A_ZY_2
\]
so
\begin{eqnarray*}
[S_XA_XY,Z]&=&[S_ZA_ZY_1,Z]+[A_ZY_2,Z]=[S_ZA_ZY_1,Z]+\mc{E}
\end{eqnarray*}
So we just have to prove that $[S_ZA_ZY_1,Z]\in \mc{E}$:
\begin{eqnarray*}
[S_ZA_ZY_1,Z]&=&S_Z^2A_ZY_1+\nabla_Z^vS_ZA_ZY_1\\
&=&S_Z^2A_ZY_1+S'_ZA_ZY_1+S_{Z'}A_ZY_1+\\
&&+S_ZA'_ZY_1+S_ZA_{Z'}Y_1+S_ZA_ZY_1'
\end{eqnarray*}
Where we denote $X'=\nabla_Z^hX$. Now using the formulas in (\ref{eq1}) and (\ref{eq2}) we obtain
\begin{eqnarray*}
[S_ZA_ZY_1,Z]&=&4S_Z^2A_ZY_1+A_ZY_1-A_ZA_Z^*A_ZY_1\\
&&+S_{Z'}A_ZY_1+S_ZA_{Z'}Y_1+S_ZA_ZY_1'\quad\in\mc{E}
\end{eqnarray*}

\subsection*{Proof of (4).}

Here we have to prove that $[\C,\A]\In \mc{E}$. Let $T\in\A$, and consider the Ricci equation
\begin{eqnarray*}
\scal{\nabla_T^h\nabla_U^hX-\nabla_U^h\nabla_T^hX-\nabla^h_{[T,U]}X,Y}=\scal{[S_X,S_Y]T,U}\\
\Ra \scal{-\nabla_T^hA_X^*U+\nabla_U^hA_X^*T+A_X^*[T,U],Y}=\scal{[S_X,S_Y]T,U}
\end{eqnarray*}
for basic vector fields $X,Y$. The first term vanishes since $A_X^*U=0$. For the second term, $A_X^*T=Z$ is basic, too, hence
\[
\nabla_U^hA_X^*T=\nabla_U^hZ=-A_Z^*U=0.
\]
Therefore, we obtained the formula
\[
\scal{[U,T],A_XY}=\scal{[S_X,S_Y]T,U}
\]
Now, the term on the right is constant, hence $\scal{[U,T],\A}=const$. In order to finish the proof, we need to show that $\scal{[U,T],U}$ is constant as well.
\\

Consider the orthonormal basis $\{U, V, W\}$, where $V,W$ span $\A$. In what follows we will work in $\vf(\leaf)/\E$. Equivalently, we will consider everything up to an equivalence $\equiv$, where we say that if $V_1,V_2$ are vector fields along $\leaf$, $V_1\equiv V_2$ if and only if $V_1-V_2\in\E$. By abuse of language, we will also define an equivalence $\equiv$ between functions, by saying that $f\equiv g$ if and only if $f-g$ is constant. Notice that:
\begin{itemize}
\item if $f,g,h\equiv 0$ as functions, then $fU+gV+hW\equiv 0$ as vector field.
\item if $\ol{V}\equiv 0$ as vector field, then $\scal{\ol{V},V},\scal{\ol{V},W},\scal{\ol{V},U}\equiv 0$ as functions.
\item if $\ol{V}\equiv 0$ and $X$ is basic, then $S_X\ol{V}\equiv 0$ and $A_X^*\ol{V}$ is basic.
\item The goal of this section becomes to prove that $[U,V]\equiv 0$, $[U,W]\equiv 0$.
\end{itemize}
The result will follow from a sequence of simple remarks:
\\

i) Since $[\A,\A]\In \E$, then $0\equiv \scal{[V,W],V}=\scal{\nabla_VW,V}-\scal{\nabla_WV,V}=\scal{\nabla_VW,V}$. In particular
\[
\nabla_VW=\scal{\nabla_VW,V}V+\scal{\nabla_VW,W}W+\scal{\nabla_VW,U}U\equiv \scal{\nabla_VW,U}U
\]
in the same way, $\nabla_WV\equiv\scal{\nabla_WV,U}U$.
\\

ii) Since $\scal{[\A,U],\A}\equiv 0$, we have that $\scal{\nabla_UV,W}\equiv \scal{\nabla_VU,W}$, and so on. Therefore
\begin{eqnarray*}
\scal{\nabla_UV,W}&\equiv&\scal{\nabla_VU,W}\equiv-\scal{V,\nabla_UW}\equiv -\scal{V,\nabla_WU}\equiv\\
&\equiv& \scal{\nabla_WV,U}\equiv\scal{\nabla_VW,U}=-\scal{\nabla_UV,W}
\end{eqnarray*}

and by looking at the end points of this chain of equivalences, we get that they are all $\equiv 0$. By the previous case,
\[
\nabla_VW\equiv \scal{\nabla_VW,U}U\equiv 0, \qquad \nabla_WV\equiv \scal{\nabla_WV,U}U\equiv 0
\]
\\

iii) Define functions $f=\scal{\nabla_VU,V}$, $g=\scal{\nabla_WU,W}$. Using what said in the previous point, $\nabla_VU\equiv fV$, and $\nabla_WU=gW$, $\nabla_VV=-fU$ and $\nabla_WW=-gU$. Moreover, notice that the goal of this section follows if we prove that $f\equiv g \equiv 0$.
\\

iv) Consider the Codazzi equation $(\nabla_VS)_XW=(\nabla_WS)_XV$. By developing both sides, we get
\[
\nabla_V(S_XW)-S_{\nabla_V^hX}W-S_X(\nabla_VW)=\nabla_W(S_XV)-S_{\nabla_W^hX}V-S_X(\nabla_WV)
\]
Notice that $\nabla_V^hX=-A_X^*V$ is basic, and therefore $S_{\nabla^h_VX}W\equiv 0$. In the same way $S_{\nabla^h_WX}V\equiv 0$, and the equation above becomes
\[
\nabla_V(S_XW)-\nabla_W(S_XV)\equiv S_X([V,W])\equiv 0.
\]
Moreover,
\begin{eqnarray*}
\nabla_V(S_XW)&\equiv& -\scal{S_XW, V}f U+ \scal{S_XW, U}fV\\
\nabla_W(S_XV)&\equiv& -\scal{S_XV, W}g U+ \scal{S_XV, U}gW
\end{eqnarray*}
By the Codazzi equation above we then have
\[
f\scal{S_XW,U}\equiv 0, \quad g\scal{S_XV,U}\equiv 0, \quad (f-g)\scal{S_XV,W}\equiv 0.
\]

There are now two possibilities: either $f\equiv g\equiv 0$, in which case we are done, or $\scal{S_XV,W}=\scal{S_XV,U}=\scal{S_XU, W}=0$. Suppose we are in this situation. Since $X$ was arbitrary, it follows that either $f\equiv g\equiv 0$, or $\{U,V,W\}$ are eigenvectors of every shape operator. Moreover,  by the arbitrariness of the elements $\{V, W\}$ we either have that $f\equiv g\equiv 0$, or $\{U,V,W\}$ are eigenvectors of every shape operator, with $V,W$ lying in the same eigenspace. Define $S_XU=\mu_XU$, $S_XV=\lambda_X V$, $S_XW=\lambda_XW$.

Using the Codazzi equation $(\nabla_US)_XV=(\nabla_VS)_XU$, we now get $(\mu_X-\lambda_X)f=\mu_{A_X^*V}$, and in the same way $(\mu_X-\lambda_X)g=\mu_{A_X^*W}$. Therefore $f\equiv g\equiv 0$, unless $\mu_X=\lambda_X$ and all shape operators are multiples of the identity. But this is not the generic situation: in fact, suppose that along a horizontal geodesic, $\gamma$, we can write $S_{\gamma'(t)}=\lambda_t Id$. Then using the first formula from \ref{equations!} we get
\[
A_{\gamma'}A_{\gamma'}^*=(\lambda'-\lambda^2-1)Id
\]
This implies that $A_{\gamma'}$ has rank 3, which is not the case.
\\

Summing up, given an appropriate choice of leaf, we can suppose that the shape operators are not all multiples of the identity, and therefore $\scal{\nabla_UV,U}$, $\scal{\nabla_UW,W}$ are constant. This finishes the proof of Proposition \ref{caseIII}.
\qed

% ---------------------------------------------------- CASE 4 --------------------------------------------------- CASE 4-----------------------

\bigskip

\section{3-dimensional SRF on spheres: case IV}

\bigskip

In this section we assume we are in case IV, i.e. for every regular point $p\in S^n$ and for any $x\in \H_p$, $A_x$ has rank 2.

The strategy for this case is the following: first we prove that any foliation that falls into case IV has very restrictive conditions on the singular set. Namely, the singular set consists on a single two-dimensional singular leaf, which is homeomorphic to a 2-sphere. From this condition, we conclude that the only possible sphere with such a foliation is $S^6$, and using the Homogeneity Theorem \ref{HomThm} again, we show that the only such foliation is homogeneous. From the topology of the singular set, this foliation must be the one induced by the irreducible representation of $\SU(2)$ on $\RR^7$.

If every horizontal vector has rank 2 some restrictions follow immediately:
\begin{itemize}
\item Every singular leaf has dimension 2. In fact, given a 1-dimensional leaf $\leaf_1$, a regular leaf $\leaf_0$ and a horizontal geodesic $\gamma:[0,1]\ra S^n$ from $\leaf_0$ to $\leaf_1$, there would be 2 linearly independent holonomy Jacobi fields $J_1$, $J_2$ along $\gamma$ vanishing at $t=1$. The initial vectors $J_1(0)$, $J_2(0)$ would be in the kernel of $A_{\gamma'(0)}^*$ by Lemma \ref{usefullem} and therefore this operator would have rank at most 1, which contradicts our rank assumption $\H=\H_2$.
\item The singular locus is connected. Again, if $C_1$, $C_2$ were different components and $\gamma$ was a minimal (horizontal) geodesic between the components, as before there would be a 2-dimensional family of holonomy Jacobi fields, vanishing on either component. A contradiction would then follow as in the previous case.
\item Every horizontal geodesic through the regular part connecting two points in the singular set has length $\pi$. In fact if the length of such geodesic was less then $\pi$, we would have again a 2-dimensional family of holonomy Jacobi fields vanishing on either end of the geodesic.
\end{itemize}
We can conclude that the singular locus $\Sigma$ is a connected compact submanifold of $S^n$ with a 2-dimensional regular foliation. Moreover every horizontal geodesic starting from $\Sigma$ goes back to $\Sigma$ after distance $\pi$. From Proposition \ref{connectedness} it follows that $\Sigma$ is homeomorphic to a sphere. But we know from Ghys \cite{Ghys}, Haefliger \cite{Hae} and Browder \cite{Bro},  (see also D. Lu \cite{Lu}) that the only way a sphere can have a two dimensional foliation is if $\Sigma\simeq S^2$ consists of just one leaf.

\subsection{Some restrictions on the dimension of $S^n$}
The fact that the singular strata consist of one single (compact) leaf imposes many restrictions. First of all, the infinitesimal foliation at $\Sigma$ is a regular one dimensional foliation on $S^{n-3}$. In particular $n-3$ must be odd, and $n$ must be even.

Moreover we can consider a (regular) leaf away from the focal locus of $\Sigma$, whose distance from $\Sigma$ is not $\pi/2$. If we take a point $p$ in such a leaf then
\[
S_x\neq 0,\qquad \forall x\in \H_p.
\]
In fact, suppose that $S_x=0$ for some $x$, and let $v\in \V_p$ be such that $A_x^*v=0$. Then by Lemma \ref{usefullem} the holonomy Jacobi field along $\gamma(t)=\exp_p(tv)$ meets $\Sigma$ at distance $\pi/2$, which is a contradiction. Then the map
\begin{eqnarray*}
S:\H_p&\lra&Sym^2(\V_p)\\
x&\lmt&S_x
\end{eqnarray*}
is injective as a map from an $(n-3)$-dimensional space to a $6$-dimensional one. Therefore $n-3\leq 6$ and $n\leq 9$.

Thus the only possible spheres in case IV are $S^4$, $S^6$, and $S^8$. We now show that $S^4$ and $S^8$ cannot occur either, leaving $S^6$ for the next paragraph. The case of $S^4$ is easy to rule out, since the foliation would have codimension one and therefore $A\equiv 0$, and we are not in case IV.
\\

Suppose now we are in $S^8$, and fix a regular point $p$ as before. We have $\dim H_p=5$. For each $v\in \V_p$, $A^v$ is a skew-symmetric endomorphism of $\H$, and hence $A^v$ can have rank $0,2,4$. In particular every $A^v$ has some kernel.

We will look more carefully at what possible ranks $A^v$ can have. Suppose $A^{v_0}=0$ for some $v_0$, and pick another $v_1$, with $x\in \ker A^{v_1}$. This would mean $A^*_xv_1=A^*_xv_0=0$ and this contradicts the fact that $A_x^*$ has rank $2$. Therefore, no $A^v$ has rank $0$.

Suppose now that $\rk A^{v_1}=\rk A^{v_2}=2$ for some $v_1,v_2$ linearly independent. Then $\ker A^{v_1}$ and $\ker A^{v_2}$ intersect at some nonzero vector $x$, and again we would have a vector $x$ such that $\rk A_x^*=1<2$. This again contradicts our rank assumption. Therefore $\dim\left\{v\st \rk A^v=2\right\}\leq 1$.

Set $\V_4=\left\{v\in \V_p\st \rk A^v=4\right\}$, and consider $\PP(\V_4)\In \PP(\V_p)$. From what was said so far, $\PP(\V_4)$ misses at most one point from $\PP(\V_p)$ (call it $v_0$)  and therefore has dimension 2. Consider now the space $E\In \PP(V_4)\times \H_p$ defined by
\[
E=\left\{([v],x)\st A^vx=0\right\}
\]
The first projection $\pi_1:E\ra \PP(\V_4)$ is a line bundle over $\PP(\V_4)$, and therefore $E$ has dimension $3$, and the second projection to $\pi_2:E\ra\H_p$ is not surjective. But then we can pick any $x$ in the complement of $\im(\pi_2)\cup \ker A^{v_0}$, and for such an $x$ we would have $\rk A^*_x=3\neq 2$. This contradicts the hypothesis of case IV, and therefore $S^8$ cannor have such a foliation.
\\

For $S^6$ one can actually prove that the irreducible $\SU(2)$ representation on $\RR^7$ induces a foliation on $S^6$ which satisfies the hypotheses of case IV. We want to prove that there are no other possible foliations on this sphere.

\subsection{Foliations on $S^6$ in case IV}

Suppose we have a foliation on $S^6$ that falls in case IV. At a regular point $p$, the horizontal space $\H_p$ has dimension 3. Given $v\in \V_p$ the skew-symmetric endomorphism $A^v$ has rank $0$ or $2$. As in the discussion about $S^8$, the case $\rk A^v=0$ can never occur, and therefore every $A^v$ has rank 2. We then have a bijection
\[
\bar{\phi}:\PP(\V_p)\lra \PP(\H_p)
\]
such that $x=\bar{\phi}(v)$ iff $A_x^*v=0$. Actually, the map $\bar{\phi}$ comes from a linear isomorphism $\V_p\lra \H_p$.
In fact, consider the map
\begin{eqnarray*}
\V\lra&\wedge^2\H_p&\stackrel{*}{\lra} \H_p\\
v\lmt&A^v&\\
&x\wedge y&\lmt x\times y
\end{eqnarray*}
and call $\phi:\V_p\ra \H_p$ the composition. Notice that if we identify a vector $x\wedge y\in \wedge^2H_p$ with the skey-symmetric map
\[
(x\wedge y)(z):=\scal{x,z}y-\scal{y,z}x,
\]
the second map in the composition sends $\alpha\in\wedge^2\H_p$ to a vector $*\alpha\in \ker \alpha$. In particular $*\alpha$ spans $\ker \alpha$, since every nonzero map in $\wedge^2\H_p$ has a one-dimensional kernel. Therefore
\[
A^v(\phi(v))=A^v(*A^v)=0.
\]
%So, $\bar{\phi}=\PP(\phi)$; also, by polarizing the equation $A^v\phi(v)=0$, we get that $A^v\phi(w)$ is skew-symmetric in $v,w$.
Of course, we can think of $\phi$ as taking vertical vector fields to horizontal ones. Set
\[
\mc{E}:=\phi^{-1}(\B),\quad  \varphi:=\phi|_{\mc{E}}:\mc{E}\lra \B,\quad \psi:=\varphi^{-1}:\B\lra \mc{E}
\]
where $\B$ is the vector space of basic vector fields, in a neighborhood of $p$ in $\leaf_p$. We want to prove $\mc{E}$ satisfies the hypotheses of the Homogeneity Theorem \ref{HomThm}. Namely, we will prove the following:
\begin{enumerate}
\item $\varphi^*$ sends $\B$ to $\A$. In particular $\A$ has dimension $3$.
\item $\phi^{-1}(\B)=\A$, and therefore $\A=\mc{E}$.
\item $\scal{V_1,V_2}$ is constant for any $V_1,V_2\in\mc{E}$.
\item $S_X\mc{E}\In\mc{E}$, for any $X\in \B$.
\item $\mc{E}$ is a Lie algebra.
\end{enumerate}

\subsection*{Proof of (1).}

Consider $\V_p$, $\H_p$ with the inner product induced by the one in $T_pS^6$, and the product on $\wedge^2\H_p$, given by
\[
\scal{\alpha,\beta}={1\over 2}\sum_{i=1}^3\scal{\alpha(e_i),\beta(e_i)}
\]
where $\{e_i\}$ is an o.n.b. of $\H_p$. It is easy to check that the $*$-operator is an isometry if $\wedge^2\H_p$ is endowed with this inner product.

It is also easy to check that $\varphi^*(e_1\times e_2)=A_{e_1}e_2$, where $e_1, e_2$ belong to an orthonormal basis of $\B$. In fact,
\[
\scal{\varphi^*(e_1\times e_2),V}=\scal{e_1\times e_2,\varphi(V)}=\scal{*(e_1\wedge e_2),*A^V}
\]
Under the identification of $\wedge^2\H_p$ with $\wedge^2\H_p^*$ above, $e_1\wedge e_2$ sends $e_1$ to $e_2$, $e_2$ to $-e_1$, and $e_3$ to zero. The above inner product is then
\[
\scal{e_1\wedge e_2,A^V}={1\over 2}\bigg(\scal{e_2, A^Ve_1}-\scal{e_1,A^Ve_2}\bigg)=\scal{A_{e_1}e_2,V}
\]
Putting everything together, $\scal{\varphi^*(e_1\times e_2),V}=\scal{A_{e_1}e_2,V}$ and the claim is proved.

\subsection*{Proof of (2).}

We have to check that $\scal{\varphi(A_XY),Z}$ is constant for any $X,Y,Z\in\B$. Considering an orthonormal basis $\{e_1,e_2,e_3\}$ of $\B$, it is enough to show that $\scal{\varphi(A_{e_1}e_2),e_k}$ is constant for $k=1,2,3$. But
\begin{eqnarray*}
\scal{\varphi(A_{e_1}e_2),e_k}&=&\scal{\varphi(A_{e_1}e_2),e_i\times e_j}\\
&=&\scal{A_{e_1}e_2,\varphi^*(e_i\times e_j)}\\
&=&\scal{A_{e_1}e_2,A_{e_i}e_j}
\end{eqnarray*}
Since there are only 3 elements in the basis of $\B$, one of $\{e_i,e_j\}$ is either $e_1$ or $e_2$, and in that case we know that the inner product is constant.

\subsection*{Proof of (3).}

 This is an immediate consequence of $\mc{E}=\A$.

\subsection*{Proof of (4).}

Consider the self adjoint homomorphism $\varphi\varphi^*:\B\ra\B$, and the orthonormal basis $\{e_1,e_2,e_3\}$ of $\B$ consisting of eigenvectors of $\varphi\varphi^*$. Then $V_1:=A_{e_2}e_3$, $V_2:=A_{e_3}e_1$, $V_3:=A_{e_1}e_2$ are orthogonal in $\mc{E}$. For example
\begin{eqnarray*}
\scal{V_1,V_2}&=&\scal{\varphi^*e_1, \varphi^*e_2}\\
&=&\scal{e_1,\varphi\varphi^*e_2}\\
&=&\lambda_2\scal{e_1,e_2}=0
\end{eqnarray*}
Now, we want to show that $S_XA_YZ\in\mc{E}$ for $X,Y,Z\in \B$. 
It is easy to see that it is enough to prove that $S_{e_1}V_1, S_{e_1}V_2\in\mc{E}$:

i) Say $S_{e_1}V_2=f_1V_1+f_2V_2+f_3V_3$. Then
\begin{eqnarray*}
f_2&=&{1\over \|V_1\|^2}\scal{S_{e_1}V_2,V_2}={1\over \|V_2\|^2}\scal{S_{e_1}A_{e_3}e_1,A_{e_3}e_1}\\
&=&{1\over \|V_2\|^2} \scal{S_{e_1}A_{e_1}e_3,A_{e_1}e_3}
\end{eqnarray*}
is constant, and the same holds for $f_3$. As for $f_1$, we know
\[
f_1^2\|V_1\|^2+f_2^2\|V_2\|^2+f_3^2\|V_3\|^2=\|S_{e_1}V_2\|^2
\]
and since everything except possibly for $f_1$ is constant, then $f_1$ must be constant as well, which means $\scal{S_{e_1}A_{e_1}e_3,A_{e_2}e_3}$ is constant.
ii) Set $S_{e_1}V_1=g_1V_1+g_2V_2+g_3V_3$. Then
\begin{eqnarray*}
g_2&=&{1\over \|V_2\|^2}\scal{S_{e_1}A_{e_2}e_3,A_{e_3}e_1}\\
&=&-{1\over \|V_2\|^2}\scal{S_{e_1}A_{e_1}e_3,A_{e_2}e_3}
\end{eqnarray*}
which is constant, by the previous case. Finally, we see that $g_1$ is constant as follows:
\begin{eqnarray*}
\|V_1\|^2g_1&=&\scal{S_{e_1}A_{e_2}e_3,A_{e_2}e_3}\\
&\stackrel{(*)}{=}&\scal{\nabla^v_{e_1}A_{e_2}e_3+S_{e_2}A_{e_3}e_1+S_{e_3}A_{e_1}e_2,A_{e_2}e_3}\\
&=&\scal{\nabla_{e_1}A_{e_2}e_3,A_{e_2}e_3}+\scal{S_{e_2}A_{e_3}e_1,A_{e_2}e_3}+\scal{S_{e_3}A_{e_1}e_2,A_{e_2}e_3}\\
&=&{1\over 2}e_1\left(\|A_{e_2}e_3\|^2\right)-\scal{S_{e_2}A_{e_2}e_3,A_{e_1}e_3}-\scal{S_{e_3}A_{e_3}e_2,A_{e_1}e_2}
\end{eqnarray*}
where in $(*)$ we used the formula in \cite{GW}, pag.149. The last expression is a sum of constant terms: the first is constant since $\|A_{e_2}e_3\|^2$ is a basic function, and $e_1$ is a basic vector field. The last two terms are constant because of what we said in the previous case.

\subsection*{Proof of (5).}

We use the Ricci equation. Notice first that for any $V\in \mc{E},\,X\in \B$, $\nabla^h_VX$ is basic: in fact for any other $Y$ basic,
\begin{eqnarray*}
\scal{\nabla_V^hX,Y}&=&\scal{\nabla_X^hV,Y}\\
&=&-\scal{A_X^*V,Y}=\scal{V,A_XY}
\end{eqnarray*}
and the last term is constant. In particular, we have that for any $V_1,V_2\in \mc{E}$, $X\in \B$,
\[
\nabla_{V_1}^h\nabla_{V_2}^hX=A^{V_1}A^{V_2}X\in \B,
\]
and the Ricci equation for $X,Y\in \B$, $V_1,V_2\in \mc{E}$ becomes
\[
\scal{A^{V_1}A^{V_2}X-A^{V_2}A^{V_1}X-A^{[V_1,V_2]}X,Y}=\scal{[S_X,S_Y]V_1,V_2}.
\]
We can rewrite this equation as
\[
\scal{[V_1,V_2],A_XY}=\scal{[S_X,S_Y]V_1,V_2}-\scal{[A^{V_1},A^{V_2}]X,Y}.
\]
Now, the right hand side is constant for every $X,Y$ basic, by what we said in the previous cases. And since the $A_XY$ generate all of $\mc{E}$, we obtain that $[V_1,V_2]\in \mc{E}$.
\qed

% -------------------------------------------------------- CASE 5  ------------------------------------------------------------------------------ CASE 5  -----------------------------------

\bigskip

\section{3-dimensional SRF on spheres: case V}\label{last}

\bigskip

Suppose that we are in case five, i.e. for every regular point $p\in S^{n}$ and every $x\in \H_p$, $A_x$ has rank 1.
The first thing to note is the following:
\begin{lem}
If $\H_p=\H_{1,p}$ at some point $p$, then the subspace of $\V_p$ generated by $A_xy$, $x,y\in \H_p$, is one dimensional.
\end{lem}
\proof Consider the map
\begin{eqnarray*}
\Phi:\PP(\H_p)&\lra&\PP(\V)\\
{}[x]&\lmt&[\im A_x]
\end{eqnarray*}
We need to show that the image of $\Phi$ is a point, and in order to do so we will show that the differential $\Phi_*$ vanishes everywhere. In fact, consider $[x]\in \PP(\H)$, and a curve $[x(t)]$ around $[x]$. Let $y\in \H_p$ such that $A_xy\neq 0$, and $A_{x'}y\neq 0$ for any $x'$ in a neighborhood of $x$ in $\H_p$. In particular, around $[x]\in \PP(\H_p)$ we can see $\Phi$ as
\[
\Phi\left([x(t)]\right)=[A_{x(t)}y]=[A_yx(t)]
\]
But notice that since $A_y$ has rank one, $[A_yx(t)]$ is constant, and this proves the lemma.
\qed
Now let $U$ be a unit length vector field which is given by $A_xy$. $U$ is well defined up to a sign, and we can in general define it locally.
We will now divide the discussion into two parts:
\begin{itemize}
\item We classify the SRF's in case V with $\Sigma_1=\emptyset$: it turns out that the only possibility is $S^5$, with a homogeneous foliation arising from a reducible $\SO(3)$ action.
\item we prove that there are no SRF's in case V with $\Sigma_1\neq \emptyset$.
\end{itemize}
\subsection{SRF's in case V with $\Sigma_1=\emptyset$}\label{caseV.1}
Consider a point $p$ in the regular part of $S^n$, let $u:=U_p\in \V_p$, and define $v,w$ to be an orthonormal basis of $u^{\perp}\cap\V_p$. We can now think of shape operators as 3 by 3 matrices, since the orthonormal basis $\{u,v,w\}$ will be fixed from now on. Remember that since $v,w$ are orthogonal to $u$, $A_x^*v=A_x^*w=0$ for any $x\in \H_p$. Notice that if a shape operator $S_x$ happened to be of the form $diag(a,b,b)$, then the holonomy Jacobi fields along $\exp_p tx$ starting at $p$ with initial values $v,w$, would both vanish at the same time, and therefore there would be a one dimensional leaf. Since$\Sigma_1=\emptyset$, then the image of the map
\begin{eqnarray*}
S:\H_p&\lra&Sym^2(\V_p)\\
x&\lmt&S_x
\end{eqnarray*}
is at most a $4$ dimensional space, complementary to the space generated by the matrices $diag(a,b,b)$. Also, since the zero matrix is such a matrix, the map $x\lmt S_x$ needs to be injective, and in particular $\dim \H_p\leq 4$. Therefore, $n\leq 7$. Moreover, the map
\begin{eqnarray*}
A^u:\H_p&\lra&\H_p\\
x&\lmt& A_x^*u
\end{eqnarray*}
is nonsingular. In fact, if $A_x^*u=0$, then $A_x^*$ would be identically zero, and this cannot happen in case $V$. Since $A^u$ is also skew symmetric, then necessarily the dimension of $\H_p$ needs to be even. Equivalently, $n=\dim \H_p+3$ needs to be odd, and this restricts our attention to the cases $n=3,5,7$. The case $n=3$ is trivial. We will first show that $n=7$ is not possible either, and this will leave us with $S^5$ only.
\\

So, assume we are on $S^7$. As usual, we fix a point $p$ in the regular part of the foliation. First of all, define a two form
\[
\omega(x,y):=\scal{A_xy,u}=\scal{A^ux,y}
\]
on $\H_p$, where $u$ is the vector defined before. As already explained this two form is non-degenerate. In particular, we can find a symplectic basis $x_1,x_2,x_3,x_4$ such that
\[
\omega(x_1,x_3)=\omega(x_1,x_4)=\omega(x_2,x_3)=\omega(x_2,x_4)=0, \omega(x_1,x_2)=\omega(x_3,x_4)=1.
\]
Now consider an orthonormal frame $U, V_1, V_2$. This forms an orthonormal basis $u=U_p, v_1=V_{1,p},v_2=V_{2,p}$ at $\V_p$. We already saw that $A_{X}^*V_1=A^*_XV_2=0$. Using  the Ricci equation we get
\begin{eqnarray*}
\omega(x,y)\scal{[w_1,w_2],u}&=&\scal{[S_x,S_y]w_1,w_2}, \qquad w_1,w_2\in \V_p,\,x,y\in\H_p
\end{eqnarray*}
where $[w_1,w_2]=[W_1,W_2]$ for some extensions $W_i$ of $w_i$ written as combinations of $\{U, V_1,V_2\}$ with constant coefficients. Define
\[
\alpha(w_1,w_2):=\scal{[w_1,w_2],u}
\]
and hence
\begin{equation}\label{eq!}
[S_x,S_y]=\omega(x,y)\alpha.
\end{equation}
We will show that $\alpha=0$ from which it follows that $[S_x,S_y]=0$ and therefore $\H\In \Sym^2\V_p$ would be a 4-dimensional space of commuting symmetric matrices. That is impossible since there are at most 3-dimensional subspaces of commuting matrices in $\Sym^2\V_p$.

Let us prove then, that $\alpha=0$. From equation (\ref{eq!}), it follows immediately that
\begin{equation}\label{eq!!}
[S_{x_1}, S_{x_3}]=[S_{x_1},S_{x_4}]=0
\end{equation}
and
\begin{equation}\label{eq!!!}
[S_{x_3},S_{x_4}]=0 \iff \alpha=0.
\end{equation}
Suppose that, in some basis, $S_{x_1}$ can be written as
\[
S_{x_1}=\left(\begin{array}{ccc}c_1 & 0 & 0 \\0 & c_2 & 0 \\0 & 0 & c_3\end{array}\right)
\]
for distinct $c_1,c_2,c_3$. From equation (\ref{eq!!}) it follows then that $S_{x_3},S_{x_4}$ are diagonal as well, and therefore they commute and hence (\ref{eq!!!}) implies $\alpha=0$.

The only other possibility is that $S_{x_1}$ is of the form
\[
S_{x_1}=\left(\begin{array}{ccc} c_1 & 0 & 0 \\0 & c_1 & 0 \\0 & 0 & c_2\end{array}\right)
\]
where $c_1\neq c_2$ (we ruled out already the possibility that $c_1=c_2$ since it would imply the existence of 1-dimensional leaves).
Then $S_{x_3}, S_{x_4}, [S_{x_3},S_{x_4}]$ are of the form
\[
S_{x_3}=\left(\begin{array}{c|c}\ol{S}_3 & 0 \\\hline 0 & c_3\end{array}\right)\quad S_{x_4}=\left(\begin{array}{c|c}\ol{S}_4 & 0 \\\hline 0 & c_4\end{array}\right)\quad [S_{x_3},S_{x_4}]=\left(\begin{array}{c|c}[\ol{S}_3,\ol{S}_4] & 0 \\\hline 0 & 0\end{array}\right)
\]
for some symmetric matrices $\ol{S}_3,\ol{S}_4$. From equation (\ref{eq!}), it follows that $[S_{x_1},S_{x_2}]$ is a nonzero multiple of $[S_{x_3},S_{x_4}]$. The generic form of $S_{x_2}$ is
\[
S_{x_2}=\left(\begin{array}{c|c}\ol{S}_2 & \ul{b} \\\hline \ul{b}^T & d\end{array}\right)
\]
and we can compute $[S_{x_1},S_{x_2}]$ as
\[
[S_{x_1},S_{x_2}]=\left(\begin{array}{c|c}0 & (c_1-c_2)\ul{b} \\\hline (c_2-c_1)\ul{b}^T & 0\end{array}\right)
\]
Now, the only way $[S_{x_1},S_{x_2}]$ and $[S_{x_3},S_{x_4}]$ can be nonzero multiples of each other is when both are zero. But again, by (\ref{eq!!!}) this means $\alpha=0$.

\subsection{SRF's in case V with $\Sigma_1=\emptyset$ on $S^5$}
Suppose we are dealing with a SRF in $S^5$, falling in case V and without one-dimensional leaves.

If a component of $\Sigma_2$ has dimension $k$, then the infinitesimal foliation consists of a regular $1$-dimensional foliation on $S^{4-k}$. In particular $4-k$ is odd, and the only possibility is $k=3$. Therefore the singular stratum consists of a (possibly disconnected) 3-dimensional manifold.

Moreover, the codimension is two and the foliation is infinitesimally polar. From \cite{LT} it follows that for every horizontal geodesic, the crossing number
\[
c:=c(\gamma)=\sum_{p\in \gamma}\dim\leaf^{reg}-\dim\leaf_p
\]
is constant. In our case, $c$ can be easily computed since by Lemma \ref{usefullem} it is twice the number of eigenvectors of $S_{\gamma'(0)}$ that lie in $\ker A^*_{\gamma'(0)}$. We have many restrictions for this number:
\begin{itemize}
\item $c$ must be even.
\item Since every $S_x$ has at most 3 eigenvalues, $c\leq 6$.
\item If $c=6$, then all the eigenvalues of $S_{\gamma'(0)}$ are in $\ker A^*_{\gamma'(0)}$. In other words, $A^*_{\gamma'(0)}=0$, and this is not possible in case IV.
\item If $c=2$ then every geodesic meets $\Sigma_2$ in two, antipodal points. By Theorem \ref{connectedness}, $\Sigma_2\simeq S^k$. But since $\Sigma_2$ has a 2-dimensional regular foliation, the only possibility is $\Sigma_2\simeq S^2$ which contradicts the initial observation that $\Sigma_2$ is $3$-dimensional.\end{itemize}
Therefore, the only possibility is that $c=4$. We actually claim that any horizontal geodesic meeting $\Sigma_2$ perpendicularly, meets $\Sigma_2$ again after every $\pi/2$. In fact, if $\gamma$ is such a geodesic, parametrized in such a way that $\gamma(0)\in \Sigma_2$, then $\gamma(t)$ and $\gamma(-t)$ belong to the same leaf, since $\gamma'(0)$ and $-\gamma'(0)$ belong to the same infinitesimal leaf around $\Sigma_2$.

We now prove the following:
\begin{prop}
The foliation has closed leaves, and $\Sigma_2$ is foliated by totally geodesic spheres.
\end{prop}
\proof 
We first prove that the foliation has closed leaves. Suppose the foliation does not have closed leaves. Then the closure of the foliation is a 4-dimensional SRF on $S^5$, and each regular leaf is an isoparametric hypersurface. By the theory of isoparametric hypersurfaces (see for example \cite{Cec}) we know that the there are two singular leaves $M^+,M^-$ f dimensions $4-m^+$, $4-m^-$ respectively, and the shape operator of any regular leaves has $g\in\{1,2,3,4,6\}$ distinct eigenvalues, where $g,m^+,m^-$ are related by the equation
\[
8=2(n-1)=g(m^++m^-)
\]
Since $M^+,M^-$ are closures of leaves, their dimension is at least 2, and therefore $m^+,m^-\leq 2$, $m^++m^-\leq 4$, and $g=8/(m^++m^-)=2$ or $4$.
If $g=2$, then $m^+=m^-=2$ and it is known that the foliation is $(S_1^2,S^2)\star (S_2^2,S^2)$. This foliation is not the closure of any proper subfoliation. In fact if it were so, then the infinitesimal foliation $\fol_p$ at a point in one of the two $\ol{\fol}$-singular $S^2$ would have to be a dense foliations on $\nu_p^1S^2=S^2$, which we know does not exist.
In the case $g=4$, we have $m^+=m^-=1$ and by a theorem of Cecil, Chi and Jensen (\cite{CCJ}) the only such foliation is given by the action of $S^1\times SO(3)$ on $\RR^6=\RR^3\times \RR^3$, where $SO(3)$ acts diagonally on each copy of $\RR^3$, while $S^1$ rotates them. This foliation has a singular leaf $\leaf_s$ given by $\RR\PP^2$ with its canonical metric. Now, either the original foliation contains $\leaf_s$  as a leaf, or the original foliation restricted $\leaf_s$ is dense. We will show that it neither is.  If $\leaf_s$ is a leaf of the original foliation, then it is a regular leaf, and hence compact with finite fundamental group. But then all the leaves around $\leaf_s$ are a finite cover of $\leaf_s$, and therefore the foliation has closed leaves. If $\leaf_s$ is foliated by dense leaves, this dense foliation lifts to $S^3$ via the double cover $S^3\ra \leaf_s=\RR\PP^3$, and gives a regular foliation on the round sphere with dense leaves. But by the work of Grove and Gromoll \cite{GG} there are no such foliations.
\\

We now prove that the leaves in $\Sigma_2$ are totally geodesic. Fix a regular point $p$, a horizontal vector $x$, and the splitting $\V_p=\scal{U_p}\oplus\scal{U_p}^{\perp}$, where $U$ is the unit-length local vector field $A_XY$ for some $X,Y$ orthogonal frame of $\H$. As we said, the geodesic $\gamma(t)=\exp_ptx$ meets the singular strata 4 times, and in particular there are 2 eigenvectors of $S_x$ in $\ker A^*_x=\scal{U_p}^{\perp}$. but this space is 2-dimensional, so $\scal{U_p}^{\perp}$ is totally spanned by eigenvectors of $S_x$. Therefore $U_p$ is an eigenvector of $S_x$. Since we choose $p$ and $x$ arbitrarily, the vector field $U$ is always an eigenvector of every shape operator.

Fix now a unit-length horizontal geodesic $\gamma(t)$, let $X(t):=\gamma'(t)$, $Y(t)$ a horizontally parallel unit-length vector field perpendicular to $X$, and define $T(t):=A_{X(t)}Y(t)$. $T$ is clearly parallel to the vector field $U$ defined before, but it does not have constant norm: call $f(t):=\|T(t)\|^2$. From the discussion above, there is a function $\lambda(t)$ such that $S_{X(t)}T(t)=\lambda(t)T(t)$. From the equations (see \cite{GW}, p. 44,194) :
\begin{eqnarray*}
(\nabla_X^vA)_XY&=&2S_XA_XY\\
(\nabla_X^vS)_XT&=&S^2_XT+R(T,X)X-A_XA_X^*T
\end{eqnarray*}
we get:
\begin{eqnarray*}
T':=\nabla_X^vT&=&\nabla_X^v(A_XY)=2S_X(A_XY)=2S_XT=2\lambda T\\
f'&=&X\scal{T,T}=2\scal{T',T}=2\scal{2\lambda T,T}=4\lambda\scal{T,T}=4\lambda f\\
\lambda' T&=&(\nabla_XS)_XT=S^2_XT+R(T,X)X-A_XA_X^*T=\lambda^2T+T-A_XA_X^*T
\end{eqnarray*}
As for the last term in the last equation, notice that $\scal{A_X^*T,X}=0$ and $\scal{A_X^*T,Y}=\scal{T,T}=f$, so $AX^*T=fY$ and $$A_XA_X^*T=fA_XY=fT.$$

Summing up, we get a system of nonlinear differential equations
\[
\left\{\begin{array}{l}f'=4\lambda f\\ \lambda'=\lambda^2+1-f\end{array}\right.
\]
We have two kinds of solutions. One is the constant solution $f\equiv 1$, $\lambda\equiv 0$. The others have the form
\[
\left\{\begin{array}{l}\lambda(t)={\sin 2t\over c+ \cos 2t}\\ \\ f(t)={c^2-1\over (c+\cos 2t)^2}\end{array}\right.
\]
when the initial conditions are $h(0)=0$, $f(0)=f_{min}={c-1\over c+1}$. Notice that the solutions have the following properties:
\begin{itemize}
\item they are periodic with period $\pi$.
\item If $f(0)$ is a minimum, then $f(\pi/2)$ is a maximum and there are no other critical points in between.
\item $f_{min}\cdot f_{max}=1$. In particular the lower $f_{min}$ is, the higher is $f_{max}$
\item $f$ is symmetric around every critical point.
\end{itemize}

Now consider the orbit space $\Delta:=S^n/\fol$. Since we are assuming that the leaves are closed, $\Delta$ is a metric space and by \cite{LT} is an orbifold (of dimension 2). By a result of Lytchack (\cite{Lyt}, corollary 1.7) and Alexandrino (\cite{Ale}, corollary 1.4), since the foliation is closed and $S^n$ is simply connected there are no exceptional leaves, and therefore the regular part of $\Delta$ is a manifold. Moreover, by looking at the infinitesimal foliation at a singular leaf, we see that the tangent cone at a singular point $\ol{p}_s$ of $\Delta$ is a cone over a quotient of $\nu^1_p\Sigma/\nu^1_p\fol$, under a group of isometries that preserves the foliation. In this case $(\nu_p^1\Sigma,\nu_p^1\fol)$ is just a 2-sphere foliated by parallel small circles, and all isometries different from the identity fix the equator. In this way, the tangent cone of $\Delta$ would be a right corner, and one of the sides would consist of regular leaves. But these leaves would have to be exceptional, and we know that there are none. Therefore the tangent cone at $\ol{p}_s$ is just  $Cone(\nu^1_p\Sigma/\nu^1_p\fol)=Cone([0,\pi])=\RR\times \RR^+$, which means $\ol{p}_s$ is (topologically) a boundary point. Therefore $\Delta$ is topologically a manifold with boundary, where the boundary corresponds to $\Sigma_2/\fol|_{\Sigma_2}$. Since the curvature of $\Delta$ is positive and the boundary is totally geodesic, it follows by Gauss-Bonnet that the boundary is connected, and therefore the same is true for $\Sigma_2$. Moreover, $\Delta$ is diffeomorphic to a disk.
\\

By O'Neill's formulas the curvature of the quotient is given by $k_{\ol{p}}:=1+3\|A_xy\|^2$, where $\{x,y\}$ are an orthonormal basis of $\H_p$ at some point $p$ over $\ol{p}$. Consider now a geodesic $\ol{\gamma}$ in $\Delta$ starting at the boundary and leaving the boundary perpendicularly. This corresponds to a horizontal geodesic starting at $\Sigma_2$ and leaving $\Sigma_2$ perpendicularly. By what we said before, such a geodesic meets $\Sigma_2$ again after $\pi/2$, and the leaves it meets at times $t$ and $-t$ are the same. In the quotient, it means that $\ol{\gamma}$ meets the boundary of $\Delta$ again after $\pi/2$, and $\ol{\gamma}(t)=\ol{\gamma}(-t)$. If $k(t)$ denotes the curvature of $\Delta$ at $\ol{\gamma}(t)$, we have that $k(t)=1+3f(t)$, where $f(t)$ is the function defined above. Since $\ol{\gamma}(t)=\ol{\gamma}(-t)$, then $f(t)=f(-t)$ and $f(0)$ is a critical point. But then all the other critical points are achieved after multiples of $\pi/2$, and at all these points $\ol{\gamma}$ is at the boundary of $\Delta$. This implies that all the critical points of the curvature are attained at the boundary. Moreover, if the minimum of the curvature of $\Delta$ is attained at $\ol{p}$, then for every geodesic $\gamma$ starting at $\ol{p}$, $\gamma(\pi/2)$ attains the maximum of the curvature of $\Delta$.

Let $\gamma_1$ be the geodesic leaving $\ol{p}_{min}$ perpendicularly, and let $\ol{p}_{max}:=\gamma_1(\pi/2)$.  We claim that every geodesic starting at $\ol{p}_{min}$ goes to $\ol{p}_{max}$ at time $\pi/2$. Suppose not, and say that there is a geodesic $\gamma_2$ starting at $\ol{p}_{min}$ and ending at $\ol{q}\neq \ol{p}_{max}$ at time $\pi/2$. Let $\gamma_0$ be the minimizing geodesic between $\ol{p}_{max}$ and $\ol{q}$. Since $\gamma_1$ hits the boundary of $\Delta$ perpendicularly, $\scal{-\gamma_1'(\pi/2),\gamma_0'(0)}>0$. But then, since $sec(\delta)>1$, by Toponogov's theorem the length of $\gamma_2$ would be strictly less than $\pi/2$, which is a contradiction.

Therefore all the geodesics starting at $\ol{p}_{min}$ end up at $\ol{p}_{max}$ at $\pi/2$. In particular, the orbifold cover $\wt{\Delta}$ of $\Delta$ is a 2-sphere with a rotationally symmetric metric. Let $g=h(t)d\theta^2+dt^2$ be the metric around $\ol{p}_{min}$ along any unit length geodesic $\gamma$. We know that
\begin{eqnarray*}
-{h''(t)\over h(t)}&=&sec(\gamma(t))=1+3f(t)=1+3{c^2-1\over (c+\cos 2t)^2}\\
h(0)&=&h(\pi/2)=0\\
h'(0)&=&-h'(\pi/2)=1
\end{eqnarray*}
Using the initial conditions $h(0)=0,h'(0)=1$ we obtain the solution
\[
h(t)={\sqrt{1+c}\cdot\sin(2t)\over 2\sqrt{\cos(2t)+c}}.
\]
Using now the condition $h'(\pi/2)=-1$ we obtain $1=\sqrt{c+1\over c-1}=\sqrt{f_{min}}$. But since $f_{min}f_{max}=1$, then $f_{max}=1$ as well, and $f$ is constantly $1$. Therefore, the curvature of $\Delta$ is constantly 4, and $\Delta=S^2({1\over2})/\ZZ_2$.
\\

As a consequence of this, for every point $\ol{p}$ in the boundary of $\Delta$, every geodesic starting at $\ol{p}$ ends at a point $\ol{q}$ which does not depend on the geodesic chosen.

In $S^n$ this means that every horizontal geodesic starting from a point $p\in\Sigma_2$ meets the same leaf $\leaf_{q}$ after time $\pi/2$, and vice versa. As a consequence, the normal 2-sphere $\nu_p^{\pi/2}(\leaf_{p})$ exponentiates to a totally geodesic 2-sphere, contained in $\leaf_{q}$. By dimension reasons, $\leaf_{q}$ is a totally geodesic 2-sphere, and by symmetry the same is true for $\leaf_{p}$. In particular, $\Sigma_2$ is foliated by totally geodesic spheres. \qed
The result above says in particular that $\Sigma_2$ is a (minimal) ruled submanifold of $S^5$. From the discussion in paragraph \ref{geomSigma1} and \cite{BDJ}, one can check that, up to rigid transformation, the only possibility for $\Sigma_2$ is to be the image of the map
\begin{eqnarray*}
S^2\times S^1&\lra& S^5\\
(v,\theta)&\lmt&(v\cos\theta,v\sin\theta)
\end{eqnarray*}
and the foliation on $\Sigma_2$ is given by the images of the $S^2\times\{\theta\}$, as $\theta$ varies in $S^1$.
This is isometric to the singular stratum of the homogeneous foliation induced by the diagonal action of $\SO(3)$ on $\RR^6=\RR^3\times\RR^3$,
\[
A\lmt\left(\begin{array}{c|c}A & 0 \\ \hline 0 & A\end{array}\right)
\]
Since the foliation has codimension 2, it is totally determined by its singular stratum, together with its foliation. In particular, our foliation is isometric to the homogeneous foliation described above.

\subsection{SRF's in case V with $\Sigma_1\neq \emptyset$}
Consider $S^n\In \RR^{n+1}$. As we said, given a component $C$ of $\Sigma_1$, we  can write
\[
C=S^{n}\cap H
\]
where the $H$ is a linear subspace of $\RR^{n+1}$. Take
\[
S_0:=C=S^{n}\cap H,\quad S_1:=S^n\cap H^{\perp}.
\]
Notice the following:
\begin{itemize}
\item $S^n$ is isometric to the spherical join $S_0\star S_1$.
\item The foliation in $S^n$ preserves $S_0$ and therefore preserves $S_1$.
\end{itemize}
If the foliation in $S_1$ is 2-dimensional, then the foliation on $S^n$ is the join
\[
(S^n,\fol)=(S_0,\fol|_{S_0})\star (S_1,\fol|_{S_1})
\]
It is easy to see that this foliation is homogeneous (under an isometric $\RR^3$-action), but it does not fall under case V. In fact, let $\gamma:[0,\pi/2]\ra S^n$ be a minimizing geodesic from a point $p_0\in S_0$, to a point $p_1$ in a 2-dimensional leaf in $S_1$ and which goes through 3 dimensional leaves. By \ref{usefulprop} $A_{\gamma'(t)}=0$ for all $t\in (0,\pi/2)$.
\\

The only possibility left is that there is some 3-dimensional leaf in $S_1$, so we can consider a point in such a leaf. We are now going to prove that $S_1$ falls into case V as well.

Notice that we can think of $\H_p$ as a subspace of $\RR^{n+1}$. Under this identification, there is a splitting
\[
\H_p=H\oplus (\H_p\cap T_pS_1).
\]
If $x\in H$, then $\exp_p\left({\pi\over 2}x\right)$ lies in $S_0$. Take now a vertical vector $v\in \V_p$, and consider $J(t)$ the holonomy Jacobi field satisfying $J(0)=v$. We can choose $v$ so that $J(\pi/2)\neq 0$. Then $J(\pi/2)$ spans the vertical space of the (one dimensional) leaf at $\exp_p\left({\pi\over 2}x\right)$, and therefore $J'(\pi/2)\in H$.
Now, we know that $J(t)=J(0)\cos t+J'(0)\sin t$, so
\[
J(\pi/2)=J'(0)=-S_xv-A_x^*v=-A^vx.
\]
where we used the fact that $S_x=0$ for $x\in H$. It follows from this, that the skew-symmetric endomorphism $A^v:\H_p\ra \H_p$ preserves $H$, and therefore it preserves $\H_p\cap T_pS_1$.

Now consider $S_1$. The horizontal distribution at a point $p\in S_1$ is exaclty $\H_p\cap T_pS_1$. Given a horizontal vector $x\in\H_p\cap T_pS_1$, the $A^*$-tensor of $S_1$, $\ol{A}^*_x$ is the projection of the $A_x^*$- tensor of $S^n$, onto $\H_p\cap T_pS_1$:
\[
\ol{A}^*_x=pr_{\H_p\cap T_pS_1}A^*_x
\]
But as we said before, $A^*_x$ already belongs to $\H_p\cap T_pS_1$, therefore $\ol{A}^*_x=A^*_x$, and in particular they have the same rank. Since we know that in our case $A^*_x$ has always rank 1, the same must be true for $\ol{A}^*_x$, and $S_1$ falls into case V.
\\

Suppose $(S^n,\fol)$ is a foliation in case V, in which $\Sigma_1$ has $k>0$ components. We just proved that we can find another foliation $(S_1,\fol|_{S^1})$ in case V, in which $\Sigma_1$ has $k-1$ components. Continuing, we know that there is a foliation in case V, with $\Sigma_1$ connected. We will show that this situation can never happen, and therefore the only possibility in case V is that $\Sigma_1=\emptyset$.
\\

So consider the case  $(S^n,\fol)$ and $\Sigma_1=S^m$. Then the subsphere $S^{m'}=S^{n-m-1}\In S^n$ at distance $\pi/2$ from $\Sigma _1$ is preserved by the foliation, and $(S^m,\fol|_{S^m})$ is again in case V, with $\Sigma_1=\emptyset$. From section (\ref{caseV.1}) we know that the only possibility is $(S^5,SO(3))$. In particular, regular leaves are quotients of $S^3$. Consider one such leaf, $\leaf$, with trivial holonomy (we know there exists at least one such leaf) and consider now a horizontal geodesic joining $p\in \leaf$ to some leaf $\leaf_0\In \Sigma_1$. Let $x=\gamma'(0)$. Since the holonomy of $\leaf$ is trivial, there exists a global basic field $X$ with $X_p=x$, and the map $q\ra \exp_q(\pi/2X_q)$ defines a submersion $\leaf\ra \leaf_0$. But this is a contradiction, since there are no Riemannian submersions from $S^3$ to a 1-dimensional manifold.

\addcontentsline{toc}{chapter}{Bibliography}
\bibliographystyle{amsalpha} 
\bibliography{bibdraft}
\nocite{*}
\end{document}